\input amstex
\input amsppt.sty

\input epsf
\epsfverbosetrue \input amsppt.sty
\magnification 980\vsize=21 true cm \hsize=16 true cm \voffset=1.1
true cm \pageno=1 \NoRunningHeads \TagsOnRight
\def\p{\partial}
\def\ve{\varepsilon}
\def\f{\frac}

\def\na{\nabla}
\def\la{\lambda}
\def\al{\alpha}

\def\t{\tilde}

\def\o{\omega}

\def\G{\Gamma}
\def\si{\sigma}

\def\dl{\delta}

\def\q{\qquad}

\def\ds{\displaystyle}

\def\div{\operatorname{div}}
\def\sgn{\operatorname{sgn}}
\def\supp{\operatorname{supp}}

\topmatter
\topmatter \vskip 0.3 true cm \title{\bf On the blowup and
lifespan of smooth solutions to a class of 2-D  nonlinear wave
equations with small initial data}
\endtitle
\endtopmatter
\document
\vskip 0.4 true cm \footnote""{* Li Jun and Yin Huicheng were
supported by the NSFC (No.~10931007, No.~11025105, No.~11001122), by
the Doctoral Program Foundation of the Ministry of Education of China
(No.~20090091110005), and by the DFG via the joint Sino-German project
``Analysis of PDEs and application.'' This work was done when Li Jun
and Yin Huicheng were visiting the Mathematical Institute of the
University of G\"{o}ttingen.} \footnote""{** Ingo Witt was partly
supported by the DFG via the Sino-German project ``Analysis of PDEs
and application.'' } \vskip 0.3 true cm \centerline{Li,
Jun$^{1,*}$; \qquad Witt, Ingo$^{2,**}$;\qquad Yin,
Huicheng$^{1,*}$} \vskip 0.5 true cm {1. Department of Mathematics and
IMS, Nanjing University, Nanjing 210093, P.R.~China.}\vskip 0.2 true cm
{2. Mathematical Institute, University of G\"{o}ttingen,
Bunsenstr.~3-5, D-37073 G\"{o}ttingen, Germany.} \vskip 0.4 true cm
\document

\centerline {\bf Abstract} We are concerned with a
class of two-dimensional nonlinear wave equations
$\p_t^2u-\div(c^2(u)\na u)=0$ or $\p_t^2u-c(u)\div(c(u)\na u)=0$ with small initial data $(u(0,x),
\p_tu(0,x))=(\ve u_0(x), \ve u_1(x))$, where $c(u)$ is a smooth
function, $c(0)\not =0$, $x\in\Bbb R^2$, $u_0(x), u_1(x)\in
C_0^{\infty}(\Bbb R^2)$ depend only on $r=\sqrt{x_1^2+x_2^2}$, and
$\ve>0$ is sufficiently small. Such equations arise in a
pressure-gradient model of fluid dynamics, also in a liquid crystal
model or other variational wave equations. When $c'(0)\not= 0$ or
$c'(0)=0$, $c''(0)\not= 0$, we establish blowup and determine the
lifespan of smooth solutions.

\vskip 0.3 true cm {\bf Keywords:} Nonlinear wave equation,
blowup, lifespan, Klainerman-Sobolev inequality.

\vskip 0.3 true cm {\bf 2010 Mathematical Subject Classification:}
35L65, 35J70.

\vskip 0.4 true cm \centerline{\bf \S1. Introduction and main
results.} \vskip 0.4 true cm

In this paper, we shall focus on two-dimensional nonlinear wave
equation of the form
$$
\cases
&\p_t^2u-\div(c^2(u)\na u)=0,\\
&u(0,x)=\ve u_0(x),\quad \p_tu(0,x)=\ve u_1(x),
\endcases\tag1.1
$$
where $c(u)$ is a smooth function with $c(0)\not=0$, $x\in \Bbb
R^2$, $u_0(x), u_1(x)\in C_0^{\infty}(\Bbb R^2)$ depend only on
$r=\sqrt{x_1^2+x_2^2}$, and $\ve>0$ is sufficiently small. We will
assume $c(u)=1+u+O(u^2)$ or $c(u)=1+u^2+O(u^3)$ without loss of
generality.

Equation (1.1) has an interesting physical background. In [1],
[27], a pressure-gradient model for the positive pressure function
$P$ derived from the 2-D compressible full Euler system takes the
form $\ds\p_t(\f{\p_tP}{P})-\Delta P=0$. When initial data
$(P(0,x), \p_tP(0,x))=(1+\ve P_0(x), \ve P_1(x))$ is given and one
sets $u(t,x)=\ln P(t,x)$, then one obtains $\p_t^2u-\div(e^u\na
u)=0$ with initial data $u(0,x)=\ln(1+\ve P_0(x))$ and
$\p_tu(0,x)=\ve P_1(x)(1+\ve P_0(x))^{-1}$. This is the case of
$c(u)=\exp{(\ds\f{u}{2})}$ in (1.1). By [6], [12], [22], [26], a
2-D liquid crystal equation or variational wave equation takes the
form $\p_t^2u-c(u)\div(c(u)\na u)=0$. Especially, for the nematic
liquid crystal equation, one has $c(u)=\al \cos^2u+\beta \sin^2u$
with positive constants $\al$ and $\beta$ satisfying
$\al\not=\beta$. In this case,
$c(u)=\al+(\beta-\al)\sin^2u=\al+(\beta-\al)u^2+O(u^3)$ which
essentially corresponds to $c(u)=1+u^2+O(u^3)$ in (1.1).

There has been extensive and remarkable work concerning the global
existence or blowup and lifespan of smooth solutions to
$n$-dimensional ($n\ge 2$) nonlinear wave equations of the form
$$
\cases &\dsize\sum_{i,j=0}^ng_{ij}(u, \na u)\p_{ij}^2u=f(u, \na u,
\na^2u),\\
&u(0,x)=\ve u_0(x), \quad \p_{x_0}u(0,x)=\ve u_1(x),
\endcases\tag1.2
$$
where $x_0=t$, $x=(x_1, ..., x_n)$, $\na=\na_{x_0,x}$, $g_{ij}$
and $f$ are smooth functions of their arguments that satisfy
$g_{ij}(u, \na u)=c_{ij}+O(|u|+|\na u|)$ and $f(u, \na u,
\na^2u)=O(|u|^2+|\na u|^2+|\na^2u|^2)$, respectively, the $c_{ij}$
are constants, and the linear operator
$\dsize\sum_{i,j=0}^nc_{ij}\p_{ij}^2$ is strictly hyperbolic with
respect to time $t$. For the functions $g_{ij}$ and $f$
independent of $u$, for $n\ge 4$, it has been shown that (1.2)
admits a global smooth solution (see [8], [9], [16]). For
$n=2,\,3$, the authors of [3], [5], [11], [14] obtained the global
existence if null conditions hold. Otherwise, if these null
conditions do not hold, then smooth solutions blow up in finite
time and their lifespan can explicitly be determined in terms of
the initial data (see [2], [3], [7], [10], [13], [24], and the
references therein). If the functions $g_{ij}$ and $f$ depend on
$u$ and the derivatives of $u$, then problem (1.2) is much harder
and there are some partial results on the global existence or
blowup and lifespan of smooth solution when $2\le n\le 4$ (for
$n\ge 5$, (1.2) admits a global solution, see [17], [21]). For
$2\le n\le 4$, lower bounds on the lifespan under some suitable
restrictions were obtained in [17], [18], [21], and the references
therein. We especially point out that if the equation in (1.2) has
the form $\p_t^2u-(1+u)\Delta u=0$, then, for $n=3$, the authors
of [4] and [19], [20] established the global existence of smooth
solution. These solutions, however, often exhibit a behavior at
infinity much different from that of solutions to a linear wave
equation.

In this paper, we will concentrate on the nonlinear wave equation
(1.1), show the finite-time blowup of smooth solution, and give an
explicit expression for the lifespan $T_{\ve}$ as $\ve\to 0$.  The
main result reads:

\smallskip
{\bf Theorem 1.1.} {\it Let $u_0(x)$, $u_1(x)\in C_0^{\infty}(\Bbb
R^2)$ depend only on $r=\sqrt{x_1^2+x_2^2}$. If $u_0(x)\not\equiv 0$
or $u_1(x)\not\equiv 0$, then\/ {\rm problem (1.1)} possesses a
$C^\infty$ solution for $0\leq t <T_\varepsilon$, where $T_{\ve}$
stands for the lifespan of smooth solution $u(t,x)$.}

(i) {\it For $c(u)=1+u+O(u^2)$},
$$
\lim_{\varepsilon\rightarrow 0}\varepsilon\sqrt{T_\varepsilon}
=\tau _0\equiv -\f{1}{2\min\limits_{\si}F_0'(\si)}.\tag1.3
$$

(ii) {\it For $c(u)=1+u^2+O(u^3)$},
$$
\lim_{\varepsilon\rightarrow 0}\varepsilon^2\ln T_\varepsilon =\nu
_0\equiv -\f{1}{2\min\limits_{\si}\{F_0(\si)F_0'(\si)\}}.\tag1.4
$$
{\it Here, $F_0(\si)$ is the Friedlander radiation field for to the 2-D
linear wave equation $\square w=0$ with initial data $(w(0,x),
\p_tw(0, x))=(u_0(r), u_1(r))$.}

Recall that
$F_0(\si)=\ds\f{1}{2\pi\sqrt{2}}\int_{\si}^{+\infty}\ds\f{R(s;
u_1)-R_s'(s; u_0)}{\sqrt{s-\si}}\,ds$, where $R(s;
u_i)=\int\dl(s-\langle\o,
x\rangle)u_i(x)\,dx=\int_{-\infty}^{\infty}u_i(\sqrt{s^2+y^2})\,dy$ is
the Radon transform of $u_i(r)$ ($i=0,1$).

Henceforth, we shall also assume that $u_0$, $u_1$ are supported in
the disk $B(0, M)$, where $M>0$.

\medskip

{\bf Remark 1.1.} {\it It follows from\/ {\rm Theorem 1.1} that smooth
solutions to {\rm (1.1)} blow up in finite time provided that
$u_0(x)\not\equiv 0$ or $u_1(x)\not\equiv 0$ because of $F_0\not\equiv
0$, $F_0(M)=0$, and $\ds\lim_{\si\to -\infty}F_0(\si)=0$. Further
properties of the function $F_0(\si)$ can be found in\/ {\rm [9,
Theorem~6.2.2]}.}

\medskip

{\bf Remark 1.2.} {\it For $c(u)=1+O(u^3)$, \/ {\rm (1.1)} admits a
global smooth solution\/ {\rm (}see\/ {\rm [17])}.}

\medskip

{\bf Remark 1.3.} {\it If $c(u)=1+c_1 u+O(u^2)$ or $c(u)=1+ c_2
u^2+O(u^3)$ in\/ {\rm (1.1)}, with $c_1\neq0$ and $c_2\neq0$, then one
also has finite-time blowup of smooth solutions, and one can establish
an explicit expressions for $T_{\ve}$ as in\/ {\rm Theorem 1.1}.}

\medskip

{\bf Remark 1.4.} {\it For the Cauchy problem for the 1-D liquid
crystal equation $\p_t^2u-c(u)\p_x(c(u)\p_x u)$ $=0$ with
$(u(0,x), \p_tu(0,x))=(u_0(x), u_1(x))$, it has been shown in\/ {\rm
[6, Theorem~1]} that if there exist positive constants $c_0<c_1$ and
an $m_0\in\Bbb R$ such that $c_0\le c(u)\le c_1$ for all $u\in \Bbb R$
and $c'(m_0)\not=0$, then the $C^1$-solution $u(t,x)$ with special
initial data $u_0(x)=m_0+\ve\phi(\ds\f{x}{\ve})$ and
$u_1(x)=-\operatorname{sgn}(c'(m_0))c(u_0(x))u'_0(x)$ blows up in
finite time\rom{;} here $\phi\in C_0^1(0,1)$ and $\ve>0$ is
sufficiently small. However, for small and in general smooth initial
data $(\ve u_0(x), \ve u_1(x))$ with compact support, which certainly
do not satisfy the assumptions on the initial data made in\/ {\rm [6,
Theorem~1]}, there are no blowup results and precise estimates of the
lifespan for the 1-D or 2-D liquid crystal equations available in the
literature so far. In fact, from the proof of\/ {\rm Theorem 1.1} and
{\rm Remark~2.1\,(ii)} in \S2, it follows that smooth solutions to the
2-D liquid crystal equation $\p_t^2u-c(u)\div(c(u)\na u)=0$ develop singularities in finite time for
arbitrarily given small and smooth, spherically symmetric initial data
with compact support. We shall discuss these points, including the
blowup rate as $t\to T_{\ve}-$, in a forthcoming paper.}

\medskip

{\bf Remark 1.5.} {\it From the proof of\/ {\rm Theorem 1.1}, we infer
that $\ds\lim_{t\to T_{\ve}-}\|\na
u(t, \cdot)\|_{L^{\infty}}=\infty$. Note that this is different from
the geometric blowup $\ds\lim_{t\to T_{\ve}-}\|\na^2
u(t, \cdot)\|_{L^{\infty}}=\infty$ that occurs for solutions $u\in
C^2([0,T_{\ve}]\times \Bbb R^2)$ to the nonlinear wave equation
$\dsize\sum_{i,j=0}^2g_{ij}(\na u)\p_{ij}^2u=0$ with small initial
data $(u(0,x), \p_tu(0,x))=(\ve u_0(x), \ve u_1(x))$ {\rm (}see {\rm
[2], [3])}. In this latter case, a local shock is formed at blowup
time for the unsteady potential flow equation\/ {\rm (}see {\rm [23]}
{\rm )}.}

\smallskip

Let us comment on the proof of Theorem 1.1. To show (1.3) or (1.4),
first we study the lower bound on the lifespan $T_{\ve}$ for problem
(1.1). As in [9, Chapter~6] or [4], by constructing a suitable
approximate solution $u_a$ to (1.1), and then considering the
difference of the exact solution $u$ and $u_a$, applying the
Klainerman-Sobolev inequality, and further establishing some delicate
energy estimate, we obtain the desired lower bound on the lifespan
$T_{\ve}$. On the other hand, the solution $u$ to (1.1) is spherically
symmetric for $t<T_{\ve}$ due to the spherical symmetry of the initial
data $(u_0, u_1)$. Based on this, we can change (1.1) into a $2\times
2$ equation system in the coordinates $(t, r)$, with still $u$
appearing in the coefficients. Thanks to the good properties of the
difference of the real solution $u$ and the approximate solution $u_a$
before blowup time $T_{\ve}$ we can
treat the solution $u$ and its derivatives simultaneously to obtain a
precise estimate on the upper bound of $T_{\ve}$. We point out that
the methods in this paper are partly motivated by [8] and [13], where
equations of the form $\p_t^2u-c^2(\p_tu)\Delta u=0$ with
$c'(0)\not=0$ were studied, but only an estimates of the first order
derivatives of $u$ were required.

The paper is organized as follows. In \S2, we construct an approximate
solution $u_a$ to (1.1) in the two cases $c(u)=1+u+O(u^2)$ and
$c(u)=1+u^2+O(u^3)$, and we establish some related estimates. In \S3,
we obtain the lower bound on the lifespan $T_{\ve}$ by continuous
induction studying the nonlinear equation satisfied by
$u-u_a$. In \S4, we change the second-order equation in (1.1) into a
$2\times 2$ first-order partial differential system and further
establish some delicate estimates on $u$ and $\na u$. From this,
together with the blowup lemma of H\"ormander [9, Lemma 1.3.2], we
obtain the upper bound on $T_{\ve}$ and thus complete the proof of
Theorem~1.1. Some useful auxiliary lemmas and conclusions are given in
an appendix.

In what follows, we will make use of the following notation:

$\bullet$ $Z$ stands for one of the Klainerman vector fields in the
symmetric case,
$$
\p_r, \ \p_t, \ S=t\p_t+r\p_r, \ H=r\p_t+t\p_r.
$$

$\bullet$ $\p$ represents $\p_r$ or $\p_t$.

$\bullet$ The norm $\|f\|_{L^2}$ stands for $\|f(t,\cdot)\|_{L^2(\Bbb
R^2)}$.

\pagebreak

\centerline{\bf \S2. Construction of an approximate solution to
(1.1)}\vskip 0.4 true cm

In this section, we construct a suitable approximate solution $u_a$ to
(1.1). Then the lower bound on $T_{\ve}$ is obtained, in \S3, by a
continuous induction argument through estimating the difference of the
solution $u$ and $u_{a}$. As $c(u)$ assumes the two different forms of
$c(u)=1+u+O(u^2)$ and $c(u)=1+u^2+O(u^3)$, respectively, and since the
constructions will be slightly different in these two cases, we divide
this section into two parts.

\vskip 0.3 true cm \centerline{\bf Part 1. Construction of
$u_a$ when $c(u)=1+u+O(u^2)$.}\vskip 0.3 true cm

As in [9, Chapter~6], we introduce the slow time variable
$\tau=\varepsilon\sqrt{1+t}$ and assume that the solution to (1.1)
can be approximated by
$$
\ve r^{-\f{1}{2}}V(\tau,\si),\quad r>0,
$$
where $\si=r-t$.

Let $V(\tau, \si)$ solve the equation
$$\cases
\partial_{\tau
\si}^2 V+2V\partial_{\si}^2 V+2(\partial_{\si} V)^2=0, \quad (\tau,\si)
\in \Bbb R^{+}\times \Bbb R ,\\
V(0,\si)=F_0(\si) ,\\
\supp V\subseteq \{(\tau,\si):\si\leq M\};\\
\endcases\tag2.1
$$
$F_0(\si)$ has been introduced in Theorem 1.1.

Recall that $\tau_0\equiv
-\ds\f{1}{2\min\limits_{\si}F_0'(\si)}>0$. With regard to problem
(2.1), one then has:

{\bf Lemma 2.1.} (2.1) {\it admits a $C^\infty$ solution $V(\tau,\si)$
for $0\leq\tau<\tau_0$, and $V(\tau,\si)$ blows up as
$\tau\to\tau_0-$.}

{\bf Proof.} Set $W=\partial_{\si}V$. Then it follows from (2.1) that
$$\cases
&\p_{\tau}W+2V\p_{\si}W+2W^2=0,\qquad (\tau,\si)\in \Bbb
R^{+}\times
\Bbb R,\\
&W(0,\si)=F'_0(\si).
\endcases\tag 2.2$$

The characteristic curve $\si=\si(\tau,s)$ of (2.1) emanating from
$(0,s)$ is defined by
$$\cases
&\ds\f{d\si}{d\tau}(\tau,s)=2V(\tau,\si(\tau,s)),\\
&\si(0,s)=s.
\endcases\tag 2.3$$

Along characteristic curves, one has
$$\cases
&\ds\f{dW}{d\tau}(\tau,\si(\tau,s))+2W^2(\tau,\si(\tau,s))=0,\\
&W(0,\si(0,s))=F_0'(s),\endcases$$
which yields, for $\tau<\tau_0$,
$$W(\tau,\si(\tau,s))=\f{F_0'(s)}{1+2F_0'(s)\tau}.\tag 2.4$$

Note that the equation in (2.1) is equivalent to
$\p_{\si}(\p_{\si}V+2V\p_{\si}V)=0$. Together with the boundary
condition for $V$ in (2.1), this yields
$$\cases
&\ds\f{dV}{d\tau}(\tau,\si(\tau,s))=0,\\
&V(0,\si(0,s))=F_0(s).
\endcases$$

This means
$$V(\tau,\si(\tau,s))=F_0(s).\tag 2.5$$

>From (2.5) and (2.3), one concludes
$$\si(\tau,s)=s+2F_0(s)\tau,\tag 2.6$$ which implies that
$\p_{s}\si(\tau,s)=1+2F_0'(s)\tau>0$ for $\tau<\tau_0$. Thus it
follows from the implicit function theorem that $s=s(\tau, \si)$ is a
smooth function of the variables $\tau$ and $\si$. Consequently,
$V(\tau,\si)=F_0(s(\tau,\si))$ is a smooth solution to (2.1) for
$\tau<\tau_0$, and as $\tau\to \tau_0-$, the derivative
$V_\sigma(\tau,\si)$ blows up due to (2.4). Lemma 2.1 is proved. \qed

\smallskip

>From [9, Chapter~6], one has that $F_0(\si)\in C^{\infty}(\Bbb R)$ is
supported in $(-\infty,M]$ and obeys the estimates
$$|F_0^{(k)}(\si)|\leq C_k(1+|\si|)^{-\f{1}{2}-k},\quad k\in \Bbb
N_0.\tag 2.7$$

>From (2.7), we now derive a decay estimate for $V(\tau, \si)$ in
(2.1) for $\tau<\tau_0$ and $\si\to -\infty$.

\medskip

{\bf Lemma 2.2.}  {\it For any positive constant $b<\tau_0$ and
$0\leq\tau\leq b$, the smooth solution $V$ to\/ {\rm (2.1)} satisfies
the estimates
$$|Z^{\al}\p_{\tau}^{l}\p_{\si}^{m}V(\tau,\si)|\leq C_{\al
b}^{lm}(1+|\si|)^{-\f{1}{2}-l-m},\quad \al,l,m\in \Bbb N_0,\tag2.8$$
where $C_{\al b}^{lm}$ are positive constants depending on $b, \al, l$
and $m$.}

{\bf Proof:} When $\tau\leq b$, it follows from (2.6)--(2.7) that
$\ds\f{|s|}{2}\leq |\si|\leq 2|s|$ for large $|s|$. Together with
(2.4)--(2.5), this yields
$$|V(\tau,\si)|\leq C_b(1+|\si|)^{-\f{1}{2}},\quad
|\p_{\si}V(\tau,\si)|\leq C_b(1+|\si|)^{-\f{3}{2}}.\tag 2.9$$

By (2.6) and (2.4), one has
$$\p_{\si}s(\tau,\si)=\f{1}{1+2F_0'(s)\tau}$$
and
$$\p_{\si}^2
V(\tau,\si(\tau,s))=\f{F_0''(s)}{(1+2F_0'(s)\tau)^2}-
\f{2F_0'(s)F_0''(s)}{(1+2F_0'(s)\tau)^3},$$
which yields
$$|\p_{\si}^2 V(\tau,\si)|\leq C_b(1+|\si|)^{-\f{5}{2}}.\tag
2.10$$

On the other hand, it follows from (2.1) and (2.10) that
$$|\p_{\tau\si}^2V(\tau,\si)|\leq C_b(1+|\si|)^{-\f{7}{2}}$$
and further
$$|\p_{\tau}V(\tau,\si)|\leq C_b(1+|\si|)^{-\f{5}{2}}.\tag 2.11$$

Based on (2.9)--(2.11), by an inductive argument, one arrives at
$$|\p_{\tau}^{l}\p_{\si}^{m}V(\tau,\si)|\leq
C_b^{lm}(1+|\si|)^{-\f{1}{2}-l-m},\quad l,m\in\Bbb N_0.$$

Due to $S=\si\p_{\si}+\ds\f{\ve t}{2\sqrt{1+t}}\p_{\tau}$ and
$H=-\si\p_{\si}+\ds\f{\ve r}{2\sqrt{1+t}}\p_{\tau}$ by
Lemma~A.1\,(ii), one analogously obtains
$$|Z^{\al}\p_{\tau}^{l}\p_{\si}^{m}V(\tau,\si)|\leq C_{\al b}^{lm}
(1+|\si|)^{-\f{1}{2}-l-m},\quad \al,l,m\in\Bbb N_0,$$
and this completes the proof of Lemma~2.2. \qed

\smallskip

Next, we construct an approximate solution $u_a^I$ to
(1.1) for $0\leq\tau=\ve\sqrt{1+t}<\tau_0$.

Let $w_0$ be the solution of the linear wave equation
$$\cases
\p_t^2 w_0-\triangle w_0=0,\\
w_0(0,x)=u_0(x),\\
\p_t w_0(0,x)=u_1(x).
\endcases
$$
It follows from [9, Theorem~6.2.1] that, for any constants $l>0$ and
$0<m<1$,
$$\align
&|Z^{\al}(w_0(t,x)-r^{-\f12}F_0(\si))|\leq C_{\al
l}(1+t)^{-\f{3}{2}}(1+|\si|)^{\f{1}{2}},\quad r\geq lt, \tag 2.12\\
&|\p^{k}w_0(t,x)|\leq C_{km}(1+t)^{-1-|k|},\quad r\le mt.\tag 2.13
\endalign$$

Choose a $C^\infty$-function $\chi(s)$ such that $\chi(s)=1$ for
$s\leq 1$ and $\chi(s)=0$ for $s\geq 2$. For $0\leq
\tau=\ve\sqrt{1+t}<\tau_0$, we take the approximate solution
$u_a^I$ to (1.1) to be
$$
u_a^I(t,x)=\varepsilon\biggl(\chi(\varepsilon
t)w_0(t,x)+r^{-\f{1}{2}}(1-\chi(\varepsilon
t))\chi(-3\ve\si)V(\si,\tau)\biggr).\tag 2.14
$$

By Lemma 2.2 and [9, Theorem 6.2.1], one has that, for a fixed
positive constant $b<\tau_0$,
$$|Z^\al
u_a^I(t,x)|\leq C_{\al b}\ve
(1+t)^{-\f{1}{2}}(1+|\si|)^{-\f{1}{2}},\quad \tau\leq
b.\tag 2.15$$

Set $J_a^I=\p_t^2 u_a^I-c^2(u_a^I)\triangle
u_a^I-2c(u_a^I)c'(u_a^I)|\nabla u_a^I|^2$.

\medskip

{\bf Lemma 2.3.} {\it One has
$$\ds \int_0^{\f{b^2}{\ve^2}-1}\|Z^\alpha
J_a^I\|_{L^2}\,dt\leq C_{\al b}\ve^{\f{3}{2}}.
$$}

{\bf Proof.} We divide the proof into three parts.

{\bf (i) $0\leq t\leq \ds\f{1}{\varepsilon}$.}
In this case, $\chi(\ve t)=1$ and $u_a^I=\ve w_0$.  This yields
$$J_a^I=\ve(1-c^2(\ve w_0))\Delta w_0-2\ve^2 c(\ve w_0)c'(\ve w_0)|\na
w_0|^2.$$

It follows from (2.15) and a direct computation that
$$
\|Z^{\al}J_a^I\|_{L^2}\leq C\ve^2 (1+t)^{-\f{1}{2}},\quad 0\leq
t\leq \f{1}{\ve}.\tag 2.16$$

{\bf(ii)  $\ds\f{1}{\varepsilon}\leq
t\leq\ds\f{2}{\varepsilon}$.}
Now we rewrite $u_a^I$ as
$$u_a^I=\ve w_0(t,x)+\ve(1-\chi(\ve
t))\bigl(r^{-\f{1}{2}}\chi(-3\ve\si)V(\tau,\si)-w_0(t,x)\bigr).$$
Then
$$J_a^I=J_1+J_2+J_3+J_4,\tag 2.17$$
where
$$\align J_1&=(1-c^2(u_a^I))\Delta u_a^I-2c(u_a^I)c'(u_a^I)|\na
u_a^I|^2,\\ J_2&=\ve(\p_t^2-\Delta)\biggl\{(1-\chi(\ve
t))r^{-\f{1}{2}}\chi(-3\ve\si)\bigl(V(\tau,\si)-F_0(\si)\bigr)
\biggr\},\\
J_3&=\ve(\p_t^2-\Delta)\biggl\{\chi(-3\ve\si)
\bigl(r^{-\f{1}{2}}F_0(\si)-w_0(t,x)
\bigl)\biggr\},\\
J_4&=\ve(\p_t^2-\Delta)\biggl\{(1-\chi(\ve
t))(\chi(-3\ve\si)-1)w_0(t,x)\biggr\}.
\endalign$$

We treat each $J_i$ ($1\le i\le 4$) in (2.17) separately.

>From (2.15) one obtains
$$\|Z^{\al}J_1^I\|_{L^2}\leq C_{\al b}\ve^2(1+t)^{-\f{1}{2}}.\tag 2.18$$

Since
$$\align
J_2&=\ve r^{-\f{1}{2}}(\p_t-\p_r)(\p_t+\p_r)\biggl\{(1-\chi(\ve
t))\chi(-3\ve\si)\bigl(V(\tau,\si)-F_0(\si)\bigr)\biggr\}\\
& \quad -\f{\ve}{4}r^{-\f{5}{2}}(1-\chi(\ve
t))\chi(-3\ve\si)\bigl(V(\tau,\si)-F_0(\si)\bigr)\\
&=O(\ve^3)
r^{-\f{1}{2}}(1+|\si|)^{-\f{1}{2}}+O(\ve^2)
r^{-\f{1}{2}}\int_{0}^{\tau}\p_{\tau\si}^2V(s,\si)\,ds\\
&\quad +O(\ve^2)(1+t)^{-\f{1}{2}}r^{-\f{1}{2}}(1+|\si|)^{-\f{3}{2}}
+O(\ve)(1+t)^{-\f{5}{2}}
\endalign$$
and
$$
|Z^{\al}(r^{-\f{1}{2}}\int_{0}^{\tau}\p_{s\si}^2V(s,\si)ds)|
\leq C_{\al b}\ve
(1+t)^{\f{1}{2}}r^{-\f{1}{2}}(1+|\si|)^{-\f{3}{2}},
$$
one has
$$\|Z^{\al}J_2\|_{L^2}\leq C_{\al
b}\ve^2(1+t)^{-\f{1}{2}}.\tag 2.19$$

Note that $-\ds\f{2}{3\ve}\leq\si\leq M$ holds in the support of $J_3$
which implies $r\geq \ds\f{1}{3}t$. Together with (2.12), this yields
$$\align
J_3=&O(\ve^3)(1+|\si|)^{\f{1}{2}}(1+t)^{-\f{3}{2}}+O(\ve^2)\p
\bigl(r^{-\f{1}{2}}F_0(\si)-w_0\bigr)+O(\ve)(1+t)^{-\f{5}{2}}
(1+|\si|)^{-\f{1}{2}}.
\endalign$$
On the other hand, it follows from property (i) of Lemma A.1 that
$$\align
|Z^{\al}\p\bigl(r^{-\f{1}{2}}F_0(\si)-w_0(t,x)\bigr)|&\leq
C_{\al}|\p Z^{\al}\bigl(r^{-\f{1}{2}}F_0(\si)-w_0(t,x)\bigr)|\\
&\leq
C_{\al}(1+|\si|)^{-1}|ZZ^{\al}\bigl(r^{-\f{1}{2}}F_0(\si)-w_0(t,x)\bigl)|\\
&\leq C_{\al}(1+t)^{-\f{3}{2}}(1+|\si|)^{-\f{1}{2}}.
\endalign$$
One then obtains
$$\|Z^{\al}J_{3}\|_{L^2}\leq
C_{\al}\ve^2(1+t)^{-\f{1}{2}}.\tag2.20$$

Analogously, together with (2.13), one arrives at
$$\|Z^{\al}J_4\|_{L^2}\leq C_{\al b}\ve^2(1+t)^{-2}.\tag 2.21$$

Collecting (2.18)--(2.21) yields
$$\|Z^{\al}J_a^I\|_{L^2}\leq C_{\al
b}\ve^2(1+t)^{-\f{1}{2}},\quad \f{1}{\ve}\leq t\leq
\f{2}{\ve}.\tag 2.22$$

{\bf (iii)  $\ds\f{2}{\varepsilon}\leq t\leq
\f{b^2}{\ve^2}-1$}.
Together with (2.1), by a direct computation one has
$$\align
J_a^I&=-\ve^2 r^{-\f{1}{2}}\p_{\tau\si}^2\hat
V\biggl(\f{1}{\sqrt{1+t}}-r^{-\f{1}{2}}\biggr)-\ve^2
r^{-1}\biggl(\p_{\tau\si}^2\hat V+2\hat V\p_{\si}^2\hat
V+2(\p_{\si}\hat V)^2\biggr)\\
&\quad +O(\ve^3)r^{-\f{1}{2}}(1+t)^{-1}\p \hat
V+O(\ve^2)(1+t)^{-\f{3}{2}}r^{-\f{1}{2}}\p \hat V,\tag 2.23
\endalign$$
where $\hat V(\tau,\si)=\chi(-3\ve\si)V(\tau,\si)$.

It follows from (2.1) that
$$\align
\ve^2 r^{-1} & \biggl(\p_{\tau\si}^2\hat V+2\hat V\p_{\si}^2\hat
V+2(\p_{\si}\hat V)^2\biggr)\\
& = O(\ve^3)r^{-1}(1+|\si|)^{-\f{3}{2}}-\ve^2
r^{-1}\chi(-3\ve\si)(1-\chi(-3\ve\si))\p_{\tau\si}^2 V\\
&= O(\ve^3)r^{-1}(1+|\si|)^{-\f{3}{2}};\tag 2.24
\endalign$$
here we have used the fact that $\chi(-3\ve\si)(1-\chi(-3\ve\si)$ is
supported in the interval $[-\ds\f{2}{3\ve}, -\ds\f{1}{3\ve}]$.

Substituting (2.24) into (2.23) yields
$$\|Z^{\al}J_a^I\|_{L^{2}}\leq C_{\al
b}\biggl(\ve^3(1+t)^{-\f{1}{2}}+\ve^2(1+t)^{-\f{3}{2}}\biggr).\tag
2.25$$

Consequently, combining (2.16), (2.22), and (2.25), one obtains
$$
\ds\int_0^{\f{b^2}{\ve^2}-1}\|Z^\al J_a^I\|_{L^2}\,dt\leq C_{\al
b}\ve^{\f{3}{2}},
$$
which finishes the proof of Lemma 2.3. \qed

\vskip 0.3 true cm \centerline{\bf Part 2. Construction of
$u_a$ when $c(u)=1+u^2+O(u^3)$.}\vskip 0.3 true cm

When $c(u)=1+u^2+O(u^3)$, set the slow time variable to
$\tau=\ve^2\ln(1+t)$ as in [4], and assume that the solution to (1.1)
can be approximated by
$$
\ve r^{-\f{1}{2}}G(\tau,\si),\quad r>0,
$$
where $\si=r-t$, and $G(\tau, \si)$ solves the
equation
$$\cases
\partial_{\tau
\si}^2 G+G^2\partial_{\si}^2 G+2G(\partial_{\si} G)^2=0,
\quad (\tau,\si)\in \Bbb R^{+}\times \Bbb R ,\\
G(0,\si)=F_0(\si) ,\\
\supp G\subseteq \{(\tau,\si)\colon\si\leq M\}.
\endcases\tag2.26
$$

Recall that $\nu _0\equiv
-\ds\f{1}{2\min\limits_{\si}\{F_0(\si)F_0'(\si)\}}>0$. With regard
to problem (2.26), one has:

\medskip

{\bf Lemma 2.4.} {\it {\rm (2.26)} admits a $C^\infty$-solution
$G(\tau,\si)$ for $0\leq\tau<\nu_0$, and $G(\tau,\si)$ blows up as
$\tau\to\nu_0-$.}

{\bf Proof.} Set $Q=\partial_{\si}G$. Then it follows from (2.26)
that
$$\cases
&\p_{\tau}Q+G^2\p_{\si}Q+2GQ^2=0,\quad (\tau,\si)\in \Bbb
R^{+}\times
\Bbb R,\\
&Q(0,\si)=F'_0(\si),\\
&\supp Q\subseteq \{(\tau,\si)\colon\si\leq M\}.
\endcases\tag 2.27$$

The characteristic curve $\si=\si(\tau,s)$ of (2.26) emanating from
$(0,s)$ is defined by
$$\cases
&\ds\f{d\si}{d\tau}(\tau,s)=G^2(\tau,\si(\tau,s)),\\
&\si(0,s)=s.
\endcases\tag 2.28$$

Along characteristic curves, one has
$$\cases
&\ds\f{dQ}{d\tau}(\tau,\si(\tau,s))+2(GQ^2)(\tau,\si(\tau,s))=0,\\
&Q(0,\si(0,s))=F_0'(s)\endcases$$ and
$$\cases
&\ds\f{dG}{d\tau}(\tau,\si(\tau,s))=0,\\
&G(0,\si(0,s))=F_0(s).
\endcases$$

Then, for $\tau<\nu_0$,
$$\cases
&G(\tau,\si(\tau,s))=F_0(s),\\
&\ds Q(\tau,\si(\tau,s))=\f{F_0'(s)}{1+2F_0(s)F_0'(s)\tau}.
\endcases$$
Together with (2.28), this yields
$$\si(\tau,s)=s+F_0^2(s)\tau,\tag 2.29$$
which means that, for $\tau<\nu_0$,
$$\p_{s}\si(\tau,s)=1+2F_0(s)F_0'(s)\tau>0.$$ Therefore, it follows
from the implicit function theorem that $s=s(\tau, \si)$ can be taken
as a smooth function of $\tau$ and $\si$. Thus
$G(\tau,\si)=F_0(s(\tau, \si))$ is a smooth solution of (2.26) when
$\tau<\nu_0$. One then completes the proof of Lemma~2.4 as the one of
Lemma 2.1. \qed

\smallskip

Parallel to Lemma 2.2, one has:

\smallskip

{\bf Lemma 2.5.}  {\it For any positive constant $b<\nu_0$ and
$0\leq\tau\leq b$, the smooth solution $G$ of\/ {\rm (2.26)} satisfies
the estimates
$$|Z^{\al}\p_{\tau}^{l}\p_{\si}^{m}G(\tau,\si)|\leq
C_{\al b}^{lm}(1+|\si|)^{-\f{1}{2}-l-m},\quad \al,l,m\in \Bbb
N_0,\tag2.30$$}

{\bf Proof:} Since the proof is similar to that of Lemma~2.2, it is
omitted. \qed

\smallskip

Next we construct an approximate solution $u_a^{II}$ to (1.1) for
$0\leq\tau=\ve^2\ln(1+t)<\nu_0$.

As in Part 1, we take the approximate solution $u_a^{II}$ to (1.1) to
be
$$
u_a^{II}(t,x)=\varepsilon\biggl(\chi(\varepsilon
t)w_0(t,x)+r^{-\f{1}{2}}(1-\chi(\varepsilon
t))\chi(-3\ve\si)G(\tau,\si)\biggr).\tag 2.31
$$

By Lemma 2.5 and [9, Theorem 6.2.1], one has
$$|Z^\al
u_a^{II}(t,x)|\leq C_{\al b}\ve
(1+t)^{-\f{1}{2}}(1+|\si|)^{-\f{1}{2}},\quad \tau\leq b<\nu_0.\tag
2.32$$

Set $J_a^{II}=\p_t^2 u_a^{II}-c^2(u_a^{II})\triangle
u_a^{II}-2c(u_a^{II})c'(u_a^{II})|\nabla u_a^{II}|^2$. Then one has:

\medskip

{\bf Lemma 2.6.}{\it $$\ds \int_0^{ e^{\f{b}{\ve^2}-1}}\|Z^\alpha
J_a^{II}\|_{L^2}\,dt\leq C_{\al b}\ve^{\f{3}{2}}|\ln\ve|.
$$}

{\bf Proof.} We divide this proof procedure into two parts.

{\bf (i) $0\leq t\leq \ds\f{2}{\varepsilon}$.}  As in (i)--(ii) of the
proof of Lemma 2.3, one has
$$\|Z^{\al}J_a^{II}\|_{L^2}\leq C_{\al
b}\ve^2|\ln\ve|(1+t)^{-\f{1}{2}}.\tag 2.33$$ Note that the factor
$\ln\ve$ in (2.33) appears due to
$$G(\tau,\si)-F_0(\si)=\tau\int_0^{\tau}\p_sG(s, \si)\,ds
=\ve^2\ln(1+t)\int_0^{\tau}\p_sG(s,\si)\,ds=O(\ve^2|\ln \ve|), \quad
\ds\f{1}{\ve}\le t\le\ds\f{2}{\ve}.$$

{\bf (ii) $\ds\f{2}{\varepsilon}\leq t\leq \ds e^{\f{b}{\ve^2}}-1$}.
It follows from a direct computation that
$$\align
J_a^{II}&=-2\ve^3 r^{-\f{1}{2}}\p_{\tau\si}^2\hat
G\biggl(\f{1}{1+t}-\f{1}{r}\biggr)-2\ve^3
r^{-\f{3}{2}}\biggl(\p_{\tau\si}^2\hat G+{\hat G}^2\p_{\si}^2\hat
G+2\hat G(\p_{\si}\hat G)^2\biggr)\\
& \qquad +O(\ve^3)(1+t)^{-\f{5}{2}},\tag 2.34
\endalign$$
where $\hat G(\tau,\si)=\chi(-3\ve\si)G(\tau,\si)$.

Note that
$$\align
\ve^3 & r^{-\f{3}{2}}\biggl(\p_{\tau\si}^2\hat G+{\hat
G}^2\p_{\si}^2\hat
G+2\hat G(\p_{\si}\hat G)^2\biggr)\\
&= O(\ve^4)r^{-\f{3}{2}}(1+|\si|)^{-\f{3}{2}}+\ve^3
r^{-\f{3}{2}}\chi(-3\ve\si)(1-\chi^{2}(-3\ve\si))\p_{\tau\si}^2 G\\
&= O(\ve^4)(1+t)^{-1}r^{-\f{1}{2}}(1+|\si|)^{-\f{3}{2}}.
\endalign$$
One then obtains
$$\|Z^{\al}J_a^{II}\|_{L^{2}}\leq C_{\al
b}\biggl(\ve^4(1+t)^{-1}+\ve^3(1+t)^{-\f{3}{2}}\biggr).\tag 2.35$$

Combining (2.35) and (2.33) yields
$$
\ds\int_0^{e^{\f{b}{\ve^2}}-1}\|Z^\al J_a^{II}\|_{L^2}\,dt\leq
C_{\al b}\ve^{\f{3}{2}}|\ln\ve|.
$$
This completes the proof of Lemma 2.6. \qed

\medskip

{\bf Remark 2.1.} {\it Consider the 2-D variational wave equation
$\p_t^2u-c(u)\div(c(u)\na u)=0$ with initial data $(u(0,x),
\p_tu(0,x))=(\ve u_0(x), \ve u_1(x))$.

{\rm (i)} Let $c(u)=1+u+O(u^2)$. As in\/ {\rm Part 1}, $u$ can then be
approximated by $\ve r^{-\f{1}{2}}V(\tau,\si)$, where
$\tau=\ve\sqrt{1+t}$, $\si=r-t$, and $V(\tau, \si)$ solves the
equation
$$\cases
\partial_{\tau
\si}^2 V+2V\partial_{\si}^2 V+(\partial_{\si} V)^2=0, \quad (\tau,\si)
\in \Bbb R^{+}\times \Bbb R ,\\
V(0,\si)=F_0(\si) ,\\
\supp V\subseteq \{(\tau,\si)\colon \si\leq M\}.\\
\endcases\tag2.36
$$
Applying the method of characteristics, one easily proves that\/ {\rm
(2.36)} admits a smooth solution only for
$0\le\tau<\tau_0=-\ds\f{1}{\min F_0'(\si)}$.

{\rm (ii)} Let $c(u)=1+u^2+O(u^3)$. As in\/ {\rm Part 2}, $u$
can then be approximated by $\ve r^{-\f{1}{2}}G(\tau,\si)$, where $\tau=\ve^2
\ln(1+t)$, $\si=r-t$, and $G(\tau, \si)$ solves the equation
$$\cases
\partial_{\tau
\si}^2 G+G^2\partial_{\si}^2 V+G(\partial_{\si} G)^2=0, \quad (\tau,\si)\in
\Bbb R^{+}\times \Bbb R ,\\
G(0,\si)=F_0(\si) ,\\
\supp G\subseteq \{(\tau,\si)\colon \si\leq M\}.\\
\endcases\tag2.37
$$
As shown in\/ {\rm [25, Theorem~1.2]}, for some special class of
initial data $G(0,\si)$, $C^1$-solutions $G$ to\/ {\rm (2.37)} blow up
in finite time.  The Friedlander radiation field $F_0(\si)$, however,
does not meet the assumptions of\/ {\rm [25]} due to $F_0^{(k)}(M)=0$
for any $k\in\Bbb N_0$, which is different from $F_0'(M)\neq 0$
assumed in\/ {\rm [25]}. Nonetheless, one can still obtain the
finite-time blowup of smooth solution to\/ {\rm (2.37)} \rom(see\/
\rom{Lemma A.5}\rom).}

\vskip 0.5 true cm \centerline{\bf \S3. The lower bound on the
lifespan $T_\varepsilon$} \vskip 0.5 true cm

In this section, based on the preparations in \S2, we establish the
lower bound on the lifespan $T_{\ve}$ by utilizing continuous
induction and the energy method.

First we deal with the case that $c(u)=1+u+O(u^2)$ in (1.1).

\medskip

{\bf Lemma 3.1.} {\it Let $c(u)=1+u+O(u^2)$. Then, for sufficiently
small $\varepsilon$ and $0\leq \tau=\ve\sqrt{1+t}\leq b<\tau_0$,
{\rm (1.1)} admits a $C^\infty$ solution $u$ which satisfies the
estimate
$$
|Z^\kappa\p (u-u_a^{I})|\leq
C_{b}\varepsilon^{\f{3}{2}}(1+t)^{-\f12}(1+|t-r|)^{-\f{1}{2}}\tag3.1
$$
for $|\kappa|\le 2$; here $u_a^I$ has been introduced in {\rm
(2.14)}.}

{\bf Proof.} Set $v=u-u_a^{I}$. Then
$$\cases
&\p_t^2 v-c^2(u)\Delta v=F,\\
&v(0,x)=\p_t v(0,x)=0,
\endcases\tag 3.2
$$
where
$$\align
F&=-J_a^I+(c^2(u)-c^2(u_a^{I}))\Delta u_a^{I}+2c(u)c'(u)|\na v|^2
+4c(u)c'(u)\na v\cdot\na u_a^{I}\\
&\qquad +2\bigl(c(u)c'(u)-c(u_a^{I})c'(u_a^{I})\bigr)|\na
u_a^{I}|^2.\tag 3.3
\endalign
$$

We will use continuous induction to prove (3.1). To this end, we
assume that, for some $T\leq \ds\f{b^2}{\ve^2}-1$,
$$
|Z^{\kappa}\p
v|\leq\varepsilon(1+t)^{-\f12}(1+|t-r|)^{-\f{1}{2}},\quad
|\kappa|\leq 2,\enspace t\leq T, \tag 3.4
$$
holds and then prove that
$$
|Z^{\kappa}\p
v|\leq\f{1}{2}\varepsilon(1+t)^{-\f12}(1+|t-r|)^{-\f{1}{2}},\quad
|{\kappa}|\leq 2, \enspace t\leq T.\tag3.5
$$

Note that from (3.4) one has
$$|Z^{\kappa} v|\leq
C\varepsilon(1+t)^{-\f12}(1+|t-r|)^{\f{1}{2}},\quad
|{\kappa}|\leq 2, \enspace t\leq T.\tag3.6
$$

Applying  $Z^\alpha$ to both hand sides of (3.2) yields, for
$|\alpha|\leq 4$,
$$
(\p_t^2-c^2(u)\Delta)Z^\alpha v=G\equiv\sum_{|\beta|\leq
|\al|}C_{\al\beta}Z^{\beta}F+\bigl[Z^{\al},(c^2(u)-1)\Delta\bigr]v
+\sum_{|\beta|<|\al|}C_{\al\beta}'Z^{\beta}\bigl((c^2(u)-1)\Delta
v\bigr);\tag3.7
$$
here the commutator relation $[Z^{\al}, \p_t^2-\triangle]
=\dsize\sum_{|\beta|<|\al|}C''_{\al\beta}Z^{\beta}(\p_t^2-\triangle)$
with suitable constants $C_{\al\beta}, C'_{\al\beta}, C''_{\al\beta}$
has been made use of.

Next we derive an estimate of $\|\p Z^\alpha v\|_{L^2}$ from
Eq.~(3.7). Define the energy
$$
E(t)=\f{1}{2}\sum_{|\al|\leq 4}\int_{\Bbb R^2}(|\p_t Z^\alpha
v|^2+c^2(u)|\nabla Z^\al v|^2)\,dx.
$$
Multiplying both sides of (3.7) by $\p_t Z^\alpha v$ ($|\al|\le 4$),
integrating by parts in $\Bbb R^2$, and noting that $|\p u|=|\p
u_a^{I}+\p v|\leq C_b\ve(1+t)^{-\f{1}{2}}$ from the construction of
$u_{a}^{I}$ and assumption (3.4), one arrives at
$$
E'(t)\leq \f{C_b\ve}{\sqrt{1+t}}E(t)+\sum_{|\al|\leq 4}\int_{\Bbb
R^2}|G|\cdot|\p_t Z^\al v|\,dx.\tag3.8
$$
Moreover, due to the inductive hypothesis (3.4) and (2.15), one
has
$$|Z^{\kappa}u|\leq C_{b}\ve(1+t)^{-\f{1}{2}}(1+|\si|)^{\f{1}{2}}\leq C_{b}\ve,
\quad |\kappa|\leq 2, \enspace t\leq T.\tag 3.9$$

We now treat each term in the sum $\dsize\sum_{|\al|\leq 4}\int_{\Bbb
R^2}|G|\cdot|\p_t Z^\al v|\,dx$ separately.

{\bf (A) Treatment of $\dsize\sum_{|\beta|<|\al|} \int_{\Bbb
R^2}|Z^{\beta}\bigl((c^2(u)-1)\Delta v\bigr)|\cdot|\p_t Z^\alpha
v|\,dx$.}
It follows from (3.9) that, for $|\beta|<|\al|$,
$$\align
\int_{\Bbb R^2}|Z^{\beta}& ((c^2(u)-1)\Delta v)|\cdot|\p_t
Z^{\al}v|\, dx\\
&\leq C_b\sum_{|\beta_1|+|\beta_2|=|\beta|}\int_{\Bbb
R^2}|Z^{\beta_1}u|\cdot|Z^{\beta_2}\Delta v|\cdot|\p_t
Z^{\al}v|\,dx\\
&\leq C_b\sum_{|\beta_1|+|\beta_2|=|\beta|}\int_{\Bbb
R^2}|Z^{\beta_1}v|\cdot|Z^{\beta_2}\Delta v|\cdot|\p_t
Z^{\al}v|\,dx\\
& \qquad +C_b\sum_{|\beta_1|+|\beta_2|=|\beta|}\int_{\Bbb
R^2}|Z^{\beta_1}u_a^{I}|\cdot|Z^{\beta_2}\Delta v|\cdot|\p_t
Z^{\al}v|\,dx.\tag3.10
\endalign$$
The troublesome term in (3.10) is $Z^{\beta_1} v$, since a term of
the form $Z^{\beta_1} v$ might not be contained in the energy
$E(t)$. However, thanks to property (i) of Lemma A.1, one has
$$|Z^{\beta_2}\Delta v|\leq
\f{2}{1+|t-r|}\sum_{|\beta_2'|=|\beta_2|+1}|Z^{\beta_2'}\p v|.$$
Due to $|\beta|<|\al|\leq 4$, by (3.4) and Lemma A.2, the first term
in the right-hand side of (3.10) can then be estimated as
$$\align
\int_{\Bbb R^2}|Z^{\beta_1}v| & \cdot|Z^{\beta_2}\Delta v|\cdot|\p_t
Z^{\al}v|\,dx\\
& \leq C_b\sum_{|\beta_2'|=|\beta_2|+1}\int_{\Bbb
R^2}|\f{1}{1+|t-r|}Z^{\beta_1}v|\cdot|Z^{\beta'_2}\p v|\cdot|\p_t
Z^{\al}v|\,dx\\
& \leq \f{C_b\ve}{\sqrt{1+t}}E(t).\tag3.11
\endalign$$

Analogously,
$$\int_{\Bbb
R^2}|Z^{\beta_1}u_a^{I}|\cdot|Z^{\beta_2}\Delta v|\cdot|\p_t
Z^{\al}v|\,dx\le\f{C_b\ve}{\sqrt{1+t}}E(t).$$

Therefore, one obtains
$$\dsize\sum_{|\beta|<|\al|} \int_{\Bbb
R^2}|Z^{\beta}\bigl((c^2(u)-1)\Delta v\bigr)|\cdot|\p_t Z^\alpha
v|\,dx\leq \f{C_b\ve}{\sqrt{1+t}}E(t).\tag 3.12$$

{\bf (B) Treatment of $\ds\int_{\Bbb
R^2}|\bigl[Z^{\al},\bigl(c^2(u)-1\bigr)\Delta\bigr]v|\cdot|\p_t
Z^{\al}v|\,dx$.}
For $$\align \int_{\Bbb
R^2}|\bigl[Z^{\al} & ,\bigl(c^2(u)-1\bigr)\Delta\bigr]v|\cdot|\p_t
Z^{\al}v|\,dx\\
& \leq C_b\sum_{\Sb|\al_1|+|\al_2|=|\al|\\|\al_1|\geq
1\endSb}\int_{\Bbb R^2}|Z^{\al_1}u|\cdot|Z^{\al_2}\Delta
v|\cdot|\p_t Z^{\al}v|\,dx\\
& \leq C_b\biggl(\sum_{\Sb|\al_1|+|\al_2|=|\al|\\|\al_1|\geq
1\endSb}\int_{\Bbb R^2}|Z^{\al_1}u_a^{I}|\cdot|Z^{\al_2}\Delta
v|\cdot|\p_t
Z^{\al}v|\,dx\\
& \qquad +\sum_{\Sb|\al_1|+|\al_2|=|\al|\\|\al_1|\geq 1\endSb}\int_{\Bbb
R^2}|Z^{\al_1}v|\cdot|Z^{\al_2}\Delta v|\cdot|\p_t
Z^{\al}v|\,dx\biggr),
\endalign$$
by the same argument as in (3.11), one has
$$\ds\int_{\Bbb
R^2}|\bigl[Z^{\al},\bigl(c^2(u)-1\bigr)\Delta\bigr]v|\cdot|\p_t
Z^{\al}v|\,dx\leq \f{C_b\ve}{\sqrt{1+t}}E(t).\tag 3.13$$

\smallskip

Next we treat each term $\ds\int_{\Bbb R^2}|Z^{\beta}F|\cdot|\p_t
Z^{\al}v|\,dx$, $|\beta|\leq |\al|$ that is included in
$\dsize\sum_{|\al|\leq 4}\int_{\Bbb R^2}|G|\cdot|\p_t Z^\al v|\,dx$.

{\bf (C) Treatment of $\ds\int_{\Bbb
R^2}|Z^{\beta}J_a^I|\cdot|\p_t Z^{\al}v|\,dx$.}
In this case, one has
$$\int_{\Bbb R^2}|Z^{\beta}J_a^I|\cdot|\p_t Z^{\al}v|\,dx\leq
\|Z^{\beta}J_a\|_{L^{2}}\cdot\sqrt{E(t)}.\tag 3.14$$

{\bf (D) Treatment of $\ds\int_{\Bbb
R^2}|Z^{\beta}\bigl((c^2(u)-c^2(u_a^{I}))\Delta
u_a^{I}\bigr)|\cdot|\p_t Z^{\al}v|\,dx$.}
Due to (3.9) and Lemmas A.1 and A.2, a direct computation yields
$$\align
\int_{\Bbb R^2}|Z^{\beta} & \bigl((c^2(u)-c^2(u_a^{I}))\Delta
u_a^{I}\bigr)|\cdot|\p_t Z^{\al}v|\,dx\\
&\leq C_b\sum_{|\beta_1|+|\beta_2|=|\beta|}\int_{\Bbb
R^2}|Z^{\beta_1}v|\cdot|Z^{\beta_2}\Delta u_a^{I}|\cdot|\p_t
Z^{\al}v|\,dx\\
&\leq C_b \sum_{\Sb
|\beta_1|+|\beta_2|=|\beta|\\|\beta_2'|=|\beta_2|+1\endSb}\int_{\Bbb
R^2}\f{1}{1+|t-r|}|Z^{\beta_1}v|\cdot|Z^{\beta_2'}\p
u_a^{I}|\cdot|\p_t Z^{\al}v|\,dx\\
&\leq \f{C_b\ve}{\sqrt{1+t}}E(t).\tag 3.15
\endalign$$

{\bf (E) Treatment of $\ds\int_{\Bbb
R^2}|Z^{\beta}\bigl(c(u)c'(u)|\na v|^2\bigr)|\cdot|\p_t
Z^{\al}v|\,dx$.}
Similarly to (D), one has
$$\int_{\Bbb R^2}|Z^{\beta}\bigl(c(u)c'(u)|\na
v|^2\bigr)|\cdot|\p_t Z^{\al}v|dx\leq
\f{C_b\ve}{\sqrt{1+t}}E(t).\tag 3.16$$

{\bf (F) Treatment of $\ds\int_{\Bbb
R^2}|Z^{\beta}\bigl(c(u)c'(u)\na v\cdot\na
u_a^{I}\bigr)|\cdot|\p_t Z^{\al}v|\,dx$.}
It follows from a direct computation that
$$\align
\int_{\Bbb R^2}|Z^{\beta} & \bigl(c(u)c'(u)\na v\cdot\na
u_a^{I}\bigr)|\cdot|\p_t Z^{\al}v|\,dx\\
&\leq C_b\sum_{|\beta_1|+|\beta_2|\leq |\beta|}|Z^{\beta_1}\p
v|\cdot|Z^{\beta_2}\p u_a^{I}|\cdot|\p_t Z^{\al}v|\,dx\\
&\leq \f{C_b\ve}{\sqrt{1+t}}E(t).\tag 3.17
\endalign$$

{\bf (G) Treatment of $\ds\int_{\Bbb
R^2}|Z^{\beta}\bigl(\bigl(c(u)c'(u)-c(u_a^{I})c'(u_a^{I})\bigr)|\na
u_a^{I}|^2\bigr)|\cdot|\p_t Z^{\al}v|\,dx$.}
This case is also similar to (D). In particular, one has
$$\int_{\Bbb
R^2}|Z^{\beta}\bigl(\bigl(c(u)c'(u)-c(u_a^{I})c'(u_a^{I})\bigr)|\na
u_a^{I}|^2\bigr)|\cdot|\p_t Z^{\al}v|\,dx
\leq\f{C_b\ve}{\sqrt{1+t}}E(t).\tag 3.18
$$

Substituting (3.12)--(3.18) into (3.8) yields
$$
E'(t)\leq\f{C_{b}\ve}{\sqrt{1+t}}E(t)+\ds\sum_{|\beta|\leq
4}\|Z^\beta J_a^I\|_{L^2}\sqrt{E(t)}.
$$

Thus, by Lemmas 2.3 and A.3, one obtains
$$
\|\p Z^\al v\|_{L^2}\leq C_{b}\ve^{\f{3}{2}},\quad|\al|\leq 4,
$$
and further
$$
\|Z^\al\p v\|_{L^2}\leq C_{b}\ve^{\f{3}{2}},\quad|\al|\leq
4.\tag3.19
$$

By (3.19) and the Klainerman-Sobolev inequality (see [9], [16]), one
has
$$
|Z^\kappa\p v|\leq
C_{b}\ve^{\f{3}{2}}(1+t)^{-\f{1}{2}}(1+|t-r|)^{-\f{1}{2}},\quad
|\kappa|\leq 2, \enspace t\leq T,\tag 3.20
$$
which means that, for small $\ve$,
$$
|Z^\kappa\p
v|\leq\f{1}{2}\varepsilon(1+t)^{-\f{1}{2}}(1+|t-r|)^{-\f{1}{2}},\quad
|\kappa|\leq 2, \enspace t\leq T.
$$
This completes the proofs of (3.6) and (3.1). \qed

\medskip

When $c(u)=1+u^2+O(u^3)$ in (1.1), then similar to Lemma 3.1, one has:

\smallskip

{\bf Lemma 3.2.} {\it Let $c(u)=1+u^2+O(u^3)$. Then, for sufficiently
small $\ve$ and $0\leq \tau=\ve^2\ln(1+t)\leq b<\nu_0$, {\rm (1.1)}
admits a $C^{\infty}$-solution $u$ which satisfies the estimate
$$|Z^{\kappa}\p(u-u_a^{II})|\leq
C_b\ve^{\f{3}{2}}|\ln\ve|(1+t)^{-\f{1}{2}}(1+|t-r|)^{-\f{1}{2}},\tag3.21$$}
for all $|\kappa|\leq 2$.

{\bf Proof.} As in the proof of Lemma 3.1, we define the
energy
$$
E(t)=\f{1}{2}\sum_{|\al|\leq 4}\int_{\Bbb R^2}(|\p_t Z^\alpha
v|^2+c^2(u)|\nabla Z^\al v|^2)\,dx
$$
and obtain
$$E'(t)\leq \f{C_b\ve^2}{1+t}E(t)+\sum_{|\beta|\le 4}
C_b\|Z^{\beta}J_a^{II}\|_{L^2}\sqrt{E(t)}.
$$
Due to Lemmas 2.6 and A.3, one then obtains (3.21) as in the proof of
Lemma 3.1. \qed

\medskip

{\bf Remark 3.1.} Lemma 3.1 {\it implies that
$\ds\lim_{\overline{\ve\rightarrow 0}}\ve\sqrt{1+ T_\ve}\geq\tau _0$
holds for the lifespan $T_{\ve}$ of solutions to\/ {\rm (1.1)} in case
$c(u)=1+u+O(u^2)$. Hence,
$$\lim_{\overline{\ve\rightarrow 0}}\ve\sqrt{T_\ve}\geq\tau
_0.\tag3.22$$

Similarly, {\rm Lemma 3.2} implies for the lifespan $T_{\ve}$ of\/
{\rm (1.1)} in case $c(u)=1+u^2+O(u^3)$ that
$$\lim_{\overline{\ve\rightarrow 0}}\ve^2\ln T_\ve\geq\nu
_0.\tag3.23$$}

\vskip 0.3 true cm \centerline{\bf \S4. The upper bound on the
lifespan $T_\ve$}\vskip 0.3 true cm

In this section, we establish the upper bound on $T_{\ve}$. Some of
our ideas are inspired by [8] and [13]. Since, in contrast to [8],
[13], $c(u)$ in (1.1) contains the solution $u$, and not the
derivatives of $u$, our derivation has to be more careful. Thanks to
estimates of $Z^{\al}(u-u_a^I)$ and $Z^{\al}(u-u_a^{II})$ with
$|\al|\le 2$ in Lemmas 3.1 and 3.2, respectively, one observes that
$|Z^{\al}(u-u_a^I)|\le C_b\varepsilon^{\f32}(1+t)^{-\f{1}{2}}$ for
$t\le\f{b^2}{\ve^2}-1$ and $|Z^{\al}(u-u_a^{II})|\le
C_b\varepsilon^{\f32}|\ln\ve|(1+t)^{-\f{1}{2}}$ for $t\le
e^{\f{b}{\ve^2}}-1$, respectively, near the light cone. This will play a
crucial role in the analysis later on.

Set $U=r^{\f{1}{2}}u$. Because of the spherical symmetry of $u$, (1.1) can
be written as
$$\cases
\ds\partial_t^2 U-c^2(u)\p_r^2U=\f{1}{4}r^{-\f{3}{2}}c^2(u)u
+2r^{-\f{1}{2}}c(u)c'(u)\bigl(\p_r U-\f{1}{2}r^{-\f{1}{2}}u\bigr)^2, \\
\ds U(0, r)=\ve r^{\f{1}{2}}u_0,\\
\ds \p_t U(0, r)=\ve r^{\f{1}{2}}u_1.\\
\endcases\tag4.1
$$

Define the operators $L_1$ and $L_2$ by
$$
L_1=\p_t+c(u)\p_r,\q L_2=\p_t-c(u)\p_r.
$$
We also set
$$
w_1=L_2 U=(\p_t-c(u)\p_r)U,\q w_2=L_1 U=(\p_t+c(u)\p_r)U,
$$
which means $\p_t U=\ds\f{w_2+w_1}{2}$ and $\p_r
U=\ds\f{w_2-w_1}{2c(u)}$.

For
$$
L_1L_2=\p_t^2-c^2(u)\p_r^2-(L_1c(u))\p_r,\q
L_2L_1=\p_t^2-c^2(u)\p_r^2+(L_2c(u))\p_r,
$$
one has
$$
\align &L_1 w_1
=\f{1}{2r^{\f{1}{2}}c(u)}c'(u)w_1^2+\f{c'(u)}{4r^{\f{1}{2}}c(u)}
\biggl(\f{3}{r^{\f{1}{2}}}c(u)u-2w_2\biggr)w_1
+\f{1}{4r^{\f{3}{2}}}c^2(u)u+\f{1}{2r^{\f{3}{2}}}c(u)c'(u)u^2-\f{3}{4r}
c'(u)uw_2,\tag
4.2\\
&L_2
w_2=\f{1}{2r^{\f{1}{2}}c(u)}c'(u)w_2^2-\f{c'(u)}{4r^{\f{1}{2}}c(u)}
\biggl(\f{3}{r^{\f{1}{2}}}c(u)u+2w_1\biggr)w_2
+\f{1}{4r^{\f{3}{2}}}c^2(u)u+\f{1}{2r^{\f{3}{2}}}c(u)c'(u)u^2+\f{3}{4r}
c'(u)uw_1.\tag
4.3
\endalign
$$

Because of $\ds\p_r c(u)=c'(u)\p_r
u=\f{c'(u)}{2r^{\f{1}{2}}c(u)}(w_2-w_1)-\f{1}{2r}c'(u)u$, one also has
$$
\align
& L_1w_1+w_1\p_r c(u)=\f{1}{4r}c'(u)uw_1+\f{1}{4r^{\f{3}{2}}}c^2(u)u
+\f{1}{2r^{\f{3}{2}}}c(u)c'(u)u^2-\f{3}{4r}c'(u)uw_2,\\
& L_2w_2-w_2\p_r
c(u)=-\f{1}{4r}c'(u)uw_2+\f{1}{4r^{\f{3}{2}}}c^2(u)u
+\f{1}{2r^{\f{3}{2}}}c(u)c'(u)u^2+\f{3}{4r}c'(u)uw_1
\endalign
$$
and
$$
\align d(|w_1|(dr-cdt))&=\sgn w_1(L_1w_1+w_1\p_r c)\,dt\wedge dr\\
&= \sgn w_1[\f{1}{4r}c'(u)uw_1+\f{1}{4r^{\f{3}{2}}}
c^2(u)u+\f{1}{2r^{\f{3}{2}}}c(u)c'(u)u^2
-\f{3}{4r}c'(u)uw_2]\,dt\wedge dr,\tag4.4\\
d(|w_2|(dr+cdt))&= \sgn w_2(L_2w_2-w_2\p_r c)\,dt\wedge dr\\
&= \sgn w_2[-\f{1}{4r}c'(u)uw_2+\f{1}{4r^{\f{3}{2}}}
c^2(u)u+\f{1}{2r^{\f{3}{2}}}c(u)c'(u)u^2
+\f{3}{4r}c'(u)uw_1]\,dt\wedge dr.\tag4.5\\
\endalign
$$

\smallskip

We first deal with the case $c(u)=1+u+O(u^2)$.

For $c(u)=1+u+O(u^2)$, one knows from \S2 that, for
$\ve\sqrt{1+T_b}=b<\tau_0$ with $b>0$ a fixed constant, (1.1) has a
$C^\infty$-solution for $t\leq T_b$. Choose $\ve>0$ so small that
$\ds\f{1}{\ve}<\ds\f{b^2}{\ve^2}-1$. Define the characteristic curve
$\G_\la^\pm$ by $\ds\f{dr}{dt}=\pm c(u(t,r))$ and let it pass through
$(\la,0)$ in the $(r,t)$-plane. Let $D$ be the domain which is bounded
by $\G_M^+$, $\G_{\rho_0-1}^+$, $\{t=0\}$, and $\{t=T_b\}$ (see Figure
1 below), where $\rho_0$ is chosen in such a way that
$F_0'(\rho_0)=\ds\min_{\si\le M}F_0'(\si)$. Obviously, $\G_M^+$ is the
straight line $r=t+M$.

\vskip -0.25 true cm

$$\epsfysize=95mm \epsfbox{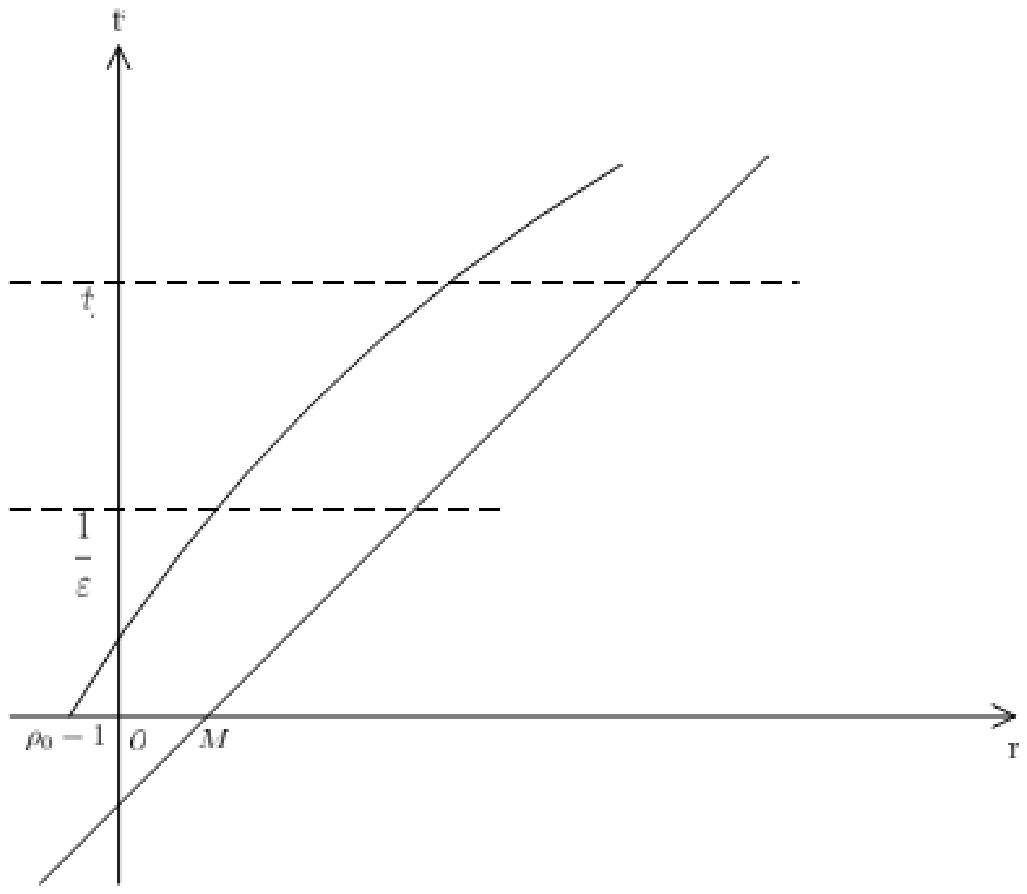}$$

\vskip -0.5 true cm
\centerline{\bf Figure 1.}

\pagebreak

One now has:

\smallskip

{\bf Lemma 4.1.} {\it If $(t,r),\,(t',r')\in\G^-_\mu\cap D$ {\rm
(}$\mu\in\Bbb R${\rm )}, and $(t,r)\in\G^+_{\la}$,
$(t',r')\in\G^+_{\la'}$, where $\la,\la'\in [\rho_0-1,M]$, then
$$
|t-t'|\leq C_b.\tag4.6
$$}

{\bf Proof.} The equation $r=r(t)$ for $\G^+_{\la}$ is
$$
\cases
&\ds\f{dr(t)}{dt}=c(u(t,r(t)))\equiv c(t),\\
&r(0)=\la,
\endcases
$$
which yields
$$r(t)-\la=\int_0^t(c(s)-1)\,ds+t.$$

Because of
$$|c(t)-1|\le C|u(t,r(t))|\le C_b\ve
(1+t)^{-\f{1}{2}}(1+|t-r(t)|)^{\f{1}{2}},\quad 0\leq
\tau=\ve \sqrt{1+t}\leq b<\tau_0,$$
one has
$$
\align
|r(t)-t|&\le |\la|+C_b\ve\int_0^t (1+s)^{-\f{1}{2}}(1+|s-r(s)|)^{\f12}\,ds\\
&\le m_0+C_b\ve\int_0^t
(1+s)^{-\f{1}{2}}(1+|s-r(s)|)^{\f12}\,ds\\
&\le m_0+C_b\ve\int_0^t (1+s)^{-\f{1}{2}}(1+|s-r(s)|)\,ds,\tag4.7
\endalign
$$
where $m_0=\max\{|\rho_0-1|,M\}$.

Set $\ds
f(t)=1+m_0+C_b\ve\int_0^{t}(1+s)^{-\f{1}{2}}(1+|s-r(s)|)ds$. By (4.7),
one then has
$$f'(t)\leq C_b\ve(1+t)^{-\f{1}{2}}f(t),\quad f(0)=1+m_0.$$
This implies, for $\ve \sqrt{1+t}\leq b$,
$$0<f(t)\le (1+m_0)e^{C_b\ve\sqrt{1+t}}\leq C_b$$
and
$$|r(t)-t|\le C_b.$$
Therefore,
$$|t+M-r(t)|\le C_b,\tag4.8$$ which means that the horizontal width
between $\G_{\rho_0-1}^+$ and $\G_M^+$ in $D$ is finite.

On the other hand, the equation $\widetilde{r}=\widetilde{r}(t)$ for
$\G^-_{\mu}$ is
$$
\cases
&\ds\f{d\widetilde{r}(t)}{dt}=-c(u(t,\widetilde{r}(t)))\equiv -\t c(t),\\
&\widetilde{r}(0)=\mu.
\endcases
$$

For $\supp u\subset \{(t,r)\colon r<t+M\}$,
$$|\t r(t)+t-\mu|\leq\cases 0,\quad t<t_{*},\\
\ds\int_{t_0}^{t}|\t c(s)+1|ds\leq
C_b\ve\int_{t_*}^{t}(1+s)^{-\f{1}{2}}(1+|\t
r(s)-s|)^{\f{1}{2}}ds,\quad t\geq t_*,\endcases$$ where $(t_*,t_*+M)$
denotes the intersection of the curves $\G^-_{\mu}$ and $\G^{+}_M$.

Combining this with (4.8) yields
$$|\t r(t)+t-\mu|\leq C_b\ve\sqrt{1+t}\leq C_b.$$

Thus, if $(t,r),(t',r')\in\G^-_\mu\cap D$ ($\mu\in\Bbb R$), and
$(t,r)\in\G^+_{\la}$, $(t',r')\in\G^+_{\la'}$, one then arrives at
$$
|t-t'|\leq\f{1}{2}\bigl(|t+r-\mu|+|t'+r'-\mu|+|t-r-\la|
+|t'-r'-\la'|+|\la-\la'|\bigr)
\leq C_b,
$$
which finishes the proof of Lemma 4.1. \qed

\smallskip
For $t\leq T_b$, define
$$
\align
& A(t)=\sup_{1/\ve\leq s\leq t}\int_{(s,r)\in D}|w_1(s,r)|\,dr,\\
& B(t)=\sup_{\Sb 1/\ve\leq s\leq t\\(s,r)\in D\endSb}s^{\f{1}{2}}|u(s,r)|,\\
& C(t)=\sup_{\Sb 1/\ve\leq s\leq t\\(s,r)\in D\endSb}s|w_2(s,r)|.
\endalign
$$

\medskip

Then one obtains:

\smallskip

{\bf Lemma 4.2.} {\it There exists a constant $E>0$ such
that, for small $\ve$,
$$A(\f{1}{\ve})\leq \ds\f{E\ve}{2}, \quad B(\f{1}{\ve})\leq E\ve, \quad
C(\f{1}{\ve})\leq E^2\ve^2.\tag4.9$$}

{\bf Proof.} Since
$w_1=r^{\f{1}{2}}\p_tu-c(u)r^{\f{1}{2}}\p_ru-\f{1}{2}r^{-\f{1}{2}}c(u)u$,
for $t\le \f{b^2}{\ve^2}-1$, one has
$$|w_1(t,r)|\le C_b\ve.\tag4.10$$
Thus, it follows from (4.8) and (4.10) that
$$\int_{(\f{1}{\ve}, r)\in D}|w_1(s,r)|\,dr\le C_b\ve.\tag4.11$$

Furthermore, because of $|u(t,r)|\le
C_b\ve(1+t)^{-\f{1}{2}}(1+|t-r|)^{\f{1}{2}}$ for $t\le
\f{b^2}{\ve^2}-1$  together with (4.8), one has, for $(\f{1}{\ve}, r)\in D$,
$$\f{1}{\ve^{\f{1}{2}}}|u(\f{1}{\ve}, r)|\le C_b\ve.\tag4.12$$

Note that
$$\align
\ds w_2(t,r)&=\ds r^{\f{1}{2}}(\p_t u+c(u)\p_r
u)+\f{1}{2}r^{-\f{1}{2}}c(u)u\\
&=r^{\f{1}{2}}\f{S+H}{t+r}u+r^{\f{1}{2}}(c(u)-1)\p_r
u+\f{1}{2}r^{-\f{1}{2}}c(u)u,\endalign$$ which implies that
$|w_2(t,r)|\le C_b\ve(1+t)^{-\f{1}{2}}$ for $(t,r)\in D$ and $\ve
\sqrt{1+t}\le b$ in view of (3.9). Together with (4.3),
this yields
$$|L_2w_2|\le\ds\f{C_b\ve^2}{1+t}.\tag4.13$$
For $w_2(t, t+M)=0$ and (4.6), from (4.13) one then obtains
$$C(\f{1}{\ve})\le C_b\ve^2.\tag4.14$$

Collecting (4.11)--(4.12) and (4.14), completes the proof of
(4.9), where $E=8(1+C_b)$. \qed

\smallskip

Based on Lemma 4.2, we will use continuous induction to estimate the
upper bound on $T_{\ve}$ when $c(u)=1+u+O(u^2)$. To this end, we
assume that, for $0\leq t\leq T'\leq T_b$,
$$
A(t)\leq E\ve,\q B(t)\leq 2E\ve,\q C(t)\leq 3E^2\ve^2.\tag 4.15
$$

We now establish:

{\bf Lemma 4.3.} {\it Under the hypothesis\/ {\rm (4.15)} and for
$\ve$ sufficiently small, one has, for $\ds\f{1}{\ve}\leq t\leq T'$,
$$
A(t)\leq \f{2}{3}E\ve,\q B(t)\leq E\ve,\q C(t)\leq
\f{5}{2}E^2\ve^2.\tag4.16
$$
}

{\bf Proof.} First we estimate $A(t)$.
By Eq.~(4.4) and Green's formula, one has, for $\ds\f{1}{\ve}\leq t\leq T'$,
$$
\align
\int_{(t,r)\in D} & |w_1(t,r)|\,dr\\
&\leq\int_{(1/\ve,r)\in D}|w_1(1/\ve,r)|\,dr\\ &\qquad
+\iint\limits_{\Sb 1/\ve\leq s\leq t\\(s,r)\in D\endSb}
|\f{1}{4r}c'(u)uw_1+\f{1}{4r^{\f{3}{2}}}c^2(u)u+\f{1}{2r^{\f{3}{2}}}c'(u)c(u)u^2
-\f{3}{4r}c'(u)uw_2|(s,r)\,dsdr\\ &\leq \f{1}{2}E\ve+\iint\limits_{\Sb
1/\ve\leq s\leq t\\ (s,r)\in
D\endSb}|\f{1}{4r}c'(u)uw_1+\f{1}{4r^{\f{3}{2}}}c^2(u)u
+\f{1}{2r^{\f{3}{2}}}c'(u)c(u)u^2
-\f{3}{4r}c'(u)uw_2|(s,r)\,dsdr.\tag4.17
\endalign
$$
By the inductive hypothesis (4.15), one has $|u(s,r)|\leq 2E\ve
s^{-\f{1}{2}}$ for $\ds\f{1}{\ve}\le s\le T'$ and $(s,r)\in D$. Note
also that $c$ is near 1 for small $\ve$, and $|r-s|\leq C_b$ holds for
$s\geq 1/\ve$. One then has $r\geq s/2$ and
$$
\align \iint\limits_{\Sb 1/\ve\leq s\leq t\\(s,r)\in
D\endSb}|\f{1}{4r}c'(u)uw_1|\,dsdr &\leq\iint\limits_{\Sb 1/\ve\leq
s\leq t\\(s,r)\in D\endSb}\f{2E\ve}{s^{\f{3}{2}}}|w_1|(s,r)\,dsdr
=2E\ve\int_{1/\ve}^t \f{1}{s^\f{3}{2}}ds\int_{(s,r)\in
D}|w_1|(s,r)\,dr\\ &\leq 4E\ve^{\f{3}{2}}A(t)\leq
4E^2\ve^{\f{5}{2}}.\tag4.18
\endalign
$$

Similarly, one has
$$\align &\iint\limits_{\Sb 1/\ve\leq s\leq t\\(s,r)\in
D\endSb}|\f{1}{4r^{\f{3}{2}}}c^2(u)u|\,dsdr \leq
2E\ve\iint\limits_{\Sb 1/\ve\leq s\leq t\\(s,r)\in
D\endSb}\f{1}{s^2}\,dsdr \leq C_bE\ve^{2},\tag4.19\\
&\iint\limits_{\Sb 1/\ve\leq s\leq t\\(s,r)\in
D\endSb}|\f{1}{2r^{\f{3}{2}}}c'(u)c(u)u^2|\,dsdr \leq C_b
E^2\ve^2\iint\limits_{\Sb 1/\ve\leq s\leq t\\(s,r)\in
D\endSb}\f{1}{s^{\f{5}{2}}}\,dsdr\leq C_b
E^2\ve^{\f{7}{2}}\tag4.20\endalign
$$
and
$$\iint\limits_{\Sb 1/\ve\leq s\leq t\\(s,r)\in
D\endSb}|\f{3}{4r}c'(u)uw_2|(s,r)\,dsdr\leq C_b
E^2\ve^2\iint\limits_{\Sb 1/\ve\leq s\leq t\\(s,r)\in
D\endSb}\f{1}{s^{\f{5}{2}}}\,dsdr\leq C_b
E^2\ve^{\f{7}{2}}.\tag4.21$$

Substituting (4.18)--(4.21) into (4.17) yields
$$
A(t)\leq\f{1}{2}E\ve+C_b E\ve^{2}+C_b E^2\ve^{\f{7}{2}},
$$
which implies $A(t)\leq\f{2}{3}E\ve$ for sufficiently small $\ve$.

\smallskip

Next we estimate $B(t)$.
Note that $u$ satisfies the equation
$$
L_2
u=L_2(r^{-\f{1}{2}}U)=r^{-\f{1}{2}}w_1+\f{1}{2r^{\f{3}{2}}}c(u)u.\tag4.22
$$
$$\epsfysize=99mm \epsfbox{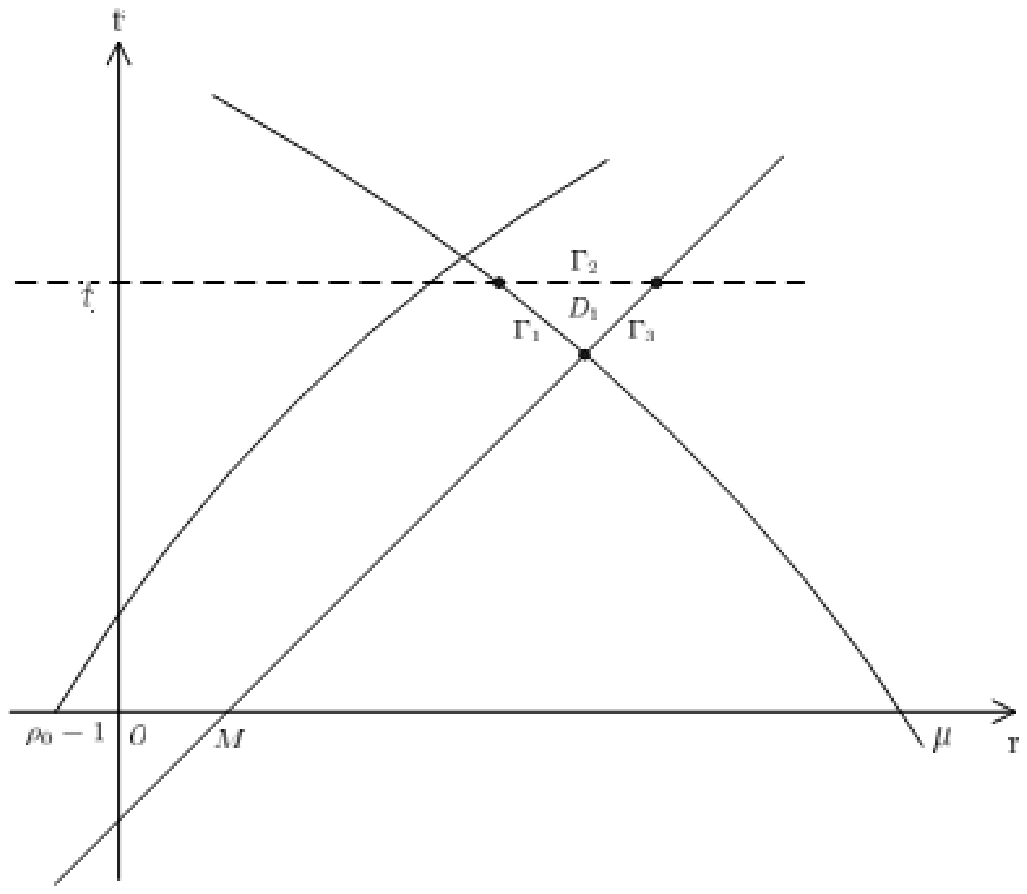}$$
\centerline{\bf Figure 2} \vskip 0.4 true cm

Integrate Eq.~(4.22) along the characteristic curve $\G_\mu^-$ which
insects $\G_M^+$ at the point $(t',r')$. When $t'\geq1/\ve$, denote by
$D_1$ the domain bounded by $\G_\mu^-$, the line $\{t=t'\}$ and
$\G_M^+$. Let $\G_1$, $\G_2$, and $\G_3$ be the faces of the boundary
(see Figure~2). By (4.4), one then has
$$
\multline
\iint\limits_{(s,r)\in D_1}\sgn
w_1[\f{1}{4r}c'(u)uw_1+\f{1}{4r^{\f{3}{2}}}c^2(u)u+\f{1}{2r^{\f{3}{2}}}
c'(u)c(u)u^2
-\f{3}{4r}c'(u)uw_2](s,r)\,dsdr \\ =(\int_{\G_1}+
\int_{\G_2})|w_1|(dr-c\,dt),
\endmultline
$$
which implies
$$
\align \int_{\G_1}|w_1|(dr-c\,ds) & \leq\iint\limits_{(s,r)\in
D_1}|\f{1}{4r}c'(u)uw_1+\f{1}{4r^{\f{3}{2}}}c^2(u)u+\f{1}{2r^{\f{3}{2}}}
c'(u)c(u)u^2 -\f{3}{4r}c'(u)uw_2|(s,r)\,dsdr\\ &\qquad
+\int_{\G_2}|w_1|(dr-c\,dt)\\ & \leq \iint\limits_{\Sb 1/\ve\leq s\leq
t\\(s,r)\in
D\endSb}|\f{1}{4r}c'(u)uw_1+\f{1}{4r^{\f{3}{2}}}c^2(u)u+\f{1}{2r^{\f{3}{2}}}
c'(u)c(u)u^2 -\f{3}{4r}c'(u)uw_2|(s,r)\,dsdr\\ &\qquad +\int_{(t,r)\in
D_1}|w_1(t,r)|\,dr\\ & \leq\f{1}{6}E\ve+\f{2}{3}E\ve=\f{5}{6}E\ve,
\endalign
$$
the last inequality resulting the estimates (4.18)--(4.21). This
yields
$$
\align \int_{t'}^t \f{|w_1(s,\t r(s))|}{(\t r(s))^{\f{1}{2}}}\,ds
&\leq\int_{t'}^t \f{|w_1(s,\t r(s))|}{(\t r(t))^{\f{1}{2}}}\,ds
\leq\f{2}{t^{\f{1}{2}}}\int_{t'}^t|w_1(s,r(s))|\,ds\\
&=\f{2}{t^{\f{1}{2}}}\int_{\G_1}\f{|w_1(s,r)|}{\sqrt{1+c^2}}\,ds
=\f{2}{t^{\f{1}{2}}}\int_{\G_1}\f{|w_1|}{2c\sqrt{1+c^2}}(dr-c\,ds)\\
&\leq\f{5}{6t^{\f{1}{2}}}E\ve.\tag 4.23
\endalign
$$

Note also that $r(s)\geq\ds\f{s}{2}\geq\ds\f{1}{2\ve}$. By (4.22), one
then has
$$
\align |u(t,r)| &\leq\int_{t'}^t \f{|w_1(s,\t r(s))|}{(\t
r(s))^{\f{1}{2}}}\,ds
+\int_{t'}^t\f{|cu|(s,\t r(s))}{2(\t r(s))^{\f{3}{2}}}\,ds\\
&\leq\f{5}{6t^{\f{1}{2}}}E\ve+3\ve^{\f{3}{2}}\int_{t'}^t
|u(s,r(s)|\,ds.
\endalign
$$

By (4.6), one obtains
$$
|u(t,\t r(t)|\leq\f{5}{6t^{\f{1}{2}}}E\ve e^{C_b\ve}\leq \ds E\ve
t^{-\f{1}{2}}.
$$

One has $t\leq t'+|t-t'|\leq 2/\ve$ for $t'\leq1/\ve$. Thus,
$t^{\f{1}{2}}|u(t,r)|\leq E\ve$ for $(t,r)\in D$, and $B(t)\le E\ve$
follows.

\smallskip

Finally, we estimate $C(t)$.
We rewrite (4.3) as
$$L_2w_2=aw_2+b, \tag4.24$$
where
$$\align
&a=\f{1}{2r^{\f{1}{2}}c(u)}c'(u)w_2-\f{c'(u)}{4r^{\f{1}{2}}c(u)}
\biggl(\f{3}{r^{\f{1}{2}}}c(u)u+2w_1\biggr),\\
&b=\f{1}{4r^{\f{3}{2}}}c^2(u)u+\f{1}{2r^{\f{3}{2}}}c(u)c'(u)u^2
+\f{3}{4}c'(u)uw_1.
\endalign
$$

Integrating (4.24) along $\G_\mu^-$ as above, one obtains
$$
|w_2(t,r)|\leq\int_{t'}^t|aw_2+b|(s,\t r(s))\,ds.\tag4.25$$

Noting $t'\geq t-|t-t'|\geq t-C_b$,
$|w_2(t,r)|\leq\ds\f{3E^2\ve^2}{t}$ and using (4.6), (4.8) and (4.23),
by the choice of $E=8(1+C_b)$ in Lemma 4.2, one arrives at
$$
\align \int_{t'}^t|b(s, & \t r(s))|\,ds\\
&\leq \int_{t'}^t\f{1}{4(\t r(s))^{\f{3}{2}}}|c^2(u)u|(s,\t
r(s))\,ds+\int_{t'}^t\f{1}{2(\t r(s))^{\f{3}{2}}}|c(u)c'(u)u^2|(s,\t r(s))\,ds\\
& \qquad +\int_{t'}^{t}\f{3}{4\t r(s)}|c'(u)uw_1|(s,\t r(s))\,ds\\
&\leq \f{1}{3}C_b
E\ve^2t^{-1}+\f{2E^2\ve^2}{t}\int_{t'}^{t}\f{1}{s^{\f{3}{2}}}\,ds+\f{5E\ve}{4t}
\int_{t'}^t \f{|w_1(s,\t r(s))|}{(\t r(s))^{\f{1}{2}}}\,ds\\
&\leq \f{1}{3}E^2\ve^2t^{-1}+\f{1}{2}E^2\ve^2
t^{-1}+\f{1}{2}E^2\ve^2 t^{-1}\\
&\leq \f{5}{3}E^2\ve^2 t^{-1},\tag4.26
\endalign
$$
and similarly
$$
\int_{t'}^t|aw_2|(s,r(s))|\,ds\leq \f{2E^2\ve^2}{t}\int_{t'}^t
|a|(s,\t r(s))\,ds \leq \f{C_b E^3\ve^3}{t}.\tag4.27
$$

Substituting (4.26)--(4.27) into (4.25) yields
$$
|w_2(t,r)|\leq\int_{t'}^t|aw_2+b|(s,\t
r(s))\,ds\leq\f{5E^2\ve^2}{2t},
$$
which shows that $C(t)\leq\ds\f{5E^2\ve^2}{2}$. \qed

\smallskip

We will use Lemma A.4 to estimate the upper bound on the lifespan
$T_{\ve}$ when $c(u)=1+u+O(u^2)$ in (1.1), based on Lemmas 4.2 and
4.3.  More specifically, we will show that
$$\overline{\lim_{\ve\rightarrow0}}\ve\sqrt{T_{\ve}}
\leq-\f{1}{2F_0'(\rho_0)}=\tau_0.\tag4.28
$$

It follows from (4.2) that on the characteristic curve
$\G_{\rho_0}^+$, $w_1(t,r(t))$ satisfies
$$\f{dw_1}{dt}(t,r(t))=L_1 w_1=a_0(t)w_1^2+a_1(t)w_1+a_2(t),\tag
4.29
$$
where
$$\align
a_0(t)&=\biggl(\f{1}{2r^{\f{1}{2}}c(u)}c'(u)\biggr)(t,r(t)),\\
a_1(t)&=\biggl(\f{c'(u)}{4r^{\f{1}{2}}c(u)}
\bigl(\f{3}{r^{\f{1}{2}}}c(u)u-2w_2\bigr)\biggr)(t,r(t)),\\
a_2(t)&=\biggl(\f{1}{4r^{\f{3}{2}}}c^2(u)u+\f{1}{2r^{\f{3}{2}}}c(u)c'(u)u^2
-\f{3}{4r}c'(u)uw_2\biggr)(t,r(t)).\endalign
$$
By (4.15), one has, for $\ds\f{1}{\ve}\leq t\leq T_b$,
$$
|a_1|\leq\f{5E\ve}{t^{\f{3}{2}}},\q|a_2|\leq\f{10E\ve^{\f{3}{2}}}{t^\f{3}{2}},
$$
which implies
$$
\int_{1/\ve}^{T_b}|a_1|\,ds\leq 10E\ve^\f{3}{2},\q
\int_{1/\ve}^{T_b}|a_2|\,ds\leq 20E^2\ve^2.\tag 4.30
$$
This also yields
$$
K=(\int_{1/\ve}^{T_b}|a_2(t)|\,dt)\exp(\int_{1/\ve}^{T_b}|a_1(t)|\,dt)
=O(\ve^2).\tag4.31
$$

By the definition of $u_a^{I}$ in (2.14), one has that
$u_a^I(1/\ve)=\ve w_0(1/\ve)$ on $\G_{\rho_0}^+$. Moreover, it follows
from Lemma 3.1 that
$$
|Z^\al(u-u_a^{I})|\leq C_{\alpha
b}\varepsilon^{\f{3}{2}}(1+t)^{-\f{1}{2}}(1+|t-r|)^{1/2}.
$$

On the other hand, by [9, Theorem 6.2.1], one has
$$
|\p^{\al}Z^{\beta}\bigl(u_a^{I}(1/\ve)-(\f{\ve}{r^{\f{1}{2}}}F_0)
(r(1/\ve)-1/\ve)\bigr)|\leq C_{\al\beta}\ve^\f{3}{2}.
$$
Therefore,
$$
\align
w_1(1/\ve) &=(r^{\f{1}{2}}\p_t
u)(1/\ve)-[c(\f{u}{2r^{\f{1}{2}}}+r^{\f{1}{2}}\p_ru)](1/\ve)\\
&=(r^{\f{1}{2}}\p_t u_a^I)(1/\ve)-[c(\f{u_a^I}{2r^{\f{1}{2}}}
+r^{\f{1}{2}}\p_ru_a^I)](1/\ve)+O(\ve^{3/2})\\ &=\ve
F_0'(r(1/\ve)-1/\ve))(-c(1/\ve)-1)+O(\ve^{3/2})\\ &=-2\ve
F_0'(r(1/\ve)-1/\ve))+O(\ve^{3/2}).
\endalign
$$

Note that one has $|r-t|\leq C+|\rho_0|$ and $|u(t,r)|\leq
C_b\ve(1+t)^{-\f{1}{2}}$ on $\G_{\rho_0}^+$. Hence,
$$
|r(1/\ve)-1/\ve-\rho_0|\leq\int_0^{1/\ve}C_b\ve(1+s)^{-\f{1}{2}}\,ds
=C_b\ve\sqrt{1+\f{1}{\ve}}.
$$

We now prove (4.28). By Lemma A.4 and (4.30)--(4.31), one has
$$
\bigl(\int_{1/\ve}^{T_b}\f{1}{2(r(t))^{\f{1}{2}}c(u)}c'(u)\,dt\bigr)
\exp(-\int_{1/\ve}^{T_b}|a_1(t)|\,dt)<(w_1(1/\ve)-K)^{-1},
$$
that is,
$$
(\sqrt{T_b}-\sqrt{1+\f{1}{\ve}})(1+O(\ve)) <(-2\ve
F_0'(\rho_0)+O(\ve^{3/2}))^{-1}\exp(5E\ve^2).
$$
Thus,
$$
\overline{\lim_{\ve\rightarrow0}}\ve\sqrt{T_{\ve}}\leq-\f{1}{2F_0'(\rho_0)}
=\tau_0,
$$
and (4.28) is shown.

\medskip

Next, we show that in case $c(u)=1+u^2+O(u^2)$ in (1.1), the lifespan
$T_{\ve}$ satisfies
$$
\overline{\lim_{\ve\rightarrow0}}\ve^2\ln T_{\ve}\leq\nu_0 =
-\f{1}{2\min\limits_{\si}\{F_0(\si)F_0'(\si)\}}.\tag4.32
$$

Although the proof is analogous to that of (4.28), for reader's
convenience, we provide the details.

Set $\t T_b=e^{\f{b}{\ve^2}}-1$, where $0<b<\nu_0$ is a fixed
constant.  As above, define $\t \G_\la^\pm$ to be the
characteristic curve given by $\ds\f{dr}{dt}=\pm c(u(t,r))$ and
passing through the point $(\la,0)$. The domain $\t D$ is bounded
by $\t \G_M^+$, $\t \G_{\t\rho_0-1}^+$, $\{t=0\}$, and $\{t=\t
T_b\}$, where $\t\rho_0$ is chosen so that
$F_0(\t\rho_0)F_0'(\t\rho_0)=\ds\min_{\si\le
M}\{F_0(\si)F_0'(\si)\}$.

\medskip

Similarly to Lemma 4.1, one has:

\smallskip

{\bf Lemma 4.4.} {\it If $(t,r), (t',r')\in \t
\Gamma_{\nu}^{-}\cap \t D$ \rom($\nu\in\Bbb R$\rom), and $(t,r)\in\t
\Gamma_{\la}^{+}, (t',r')\in\t \Gamma_{\la'}^{+}$, where
$\la,\la'\in [\t\rho_0-1,M]$, then
$$|t-t'|\leq C_b.\tag 4.33$$}

{\bf Proof.} For $\la\in [\t\rho_0-1,M]$, the equation $r=r(t)$ of
$\t \Gamma_{\la}^{+}$ is
$$\cases
&\ds\f{dr(t)}{dt}=c(u(t,r(t)))\equiv c(t),\\
&r(0)=\la.
\endcases$$

Because of $|c(t)-1|\leq C_b|u(t,r(t))|^2\leq
C_b\ve^2(1+t)^{-1}(1+|r-r(t)|)$ for $0<\tau=\ve^2\ln(1+t)\leq
b<\nu_0$, one has
$$\align
|r(t)-t|&\leq|\la|+\ds\int_{0}^{t}|c(u(s,r(s)))-1|\,ds\\
&\leq m_0+C_b\ve^2\int_0^{t}(1+s)^{-1}(1+|r(s)-s|)\,ds,
\endalign$$
which implies $|r(t)-t|\leq C_b$ for $t\le{\t T}_b$. The proof of
Lemma 4.4 then concludes by an argument similar to that in the
proof of Lemma 4.1. \qed

\medskip

Define $A(t), B(t)$ and $C(t)$ as in Lemma 4.2. When
$c(u)=1+u^2+O(u^3)$ in (1.1), one obtains paralleling Lemmas 4.2 and
4.3:

\smallskip

{\bf Lemma 4.5.} {\it There exists a positive constant $E$ such that,
for small $\ve$,

\rom{(i)} $$A(\f{1}{\ve})\leq \f{E\ve}{2},\quad B(\f{1}{\ve})\leq
E\ve,\quad C(\f{1}{\ve})\leq E^2\ve^2.\tag 4.34$$

\rom{(ii)} If $0\leq t\leq \t T_b$, then
$$
A(t)\leq E\ve,\q B(t)\leq 2E\ve,\q C(t)\leq 3E^2\ve^2.\tag 4.35
$$}

{\bf Proof.} Since the proof is analogous to those of Lemmas 4.2 and
4.3, it is omitted. \qed

\smallskip

Next we prove (4.32).  It follows from (4.2) that along the
characteristic curve $\t\G_{\t\rho_0}^+$, $w_1(t,r(t))$ satisfies
$$\f{dw_1}{dt}(t,r(t))=a_0(t)w_1^2+a_1(t)w_1+a_2(t),\tag
4.36
$$
where $a_0(t), a_1(t)$ and $a_2(t)$ are defined as in (4.29).

By (4.35), one has, for $\ds\f{1}{\ve}\leq t\leq \t T_b$,
$$
|a_1|\leq\f{5E\ve^2}{t^{2}},\q|a_2|\leq\f{10E\ve^{\f{3}{2}}}{t^\f{3}{2}},
$$
which implies
$$
\int_{1/\ve}^{\t T_b}|a_1|\,ds\leq 10E\ve^3,\q \int_{1/\ve}^{\t
T_b}|a_2|\,ds\leq 20E^2\ve^2\tag 4.37
$$
and
$$
K=(\int_{1/\ve}^{\t T_b}|a_2(t)|\,dt)\exp(\int_{1/\ve}^{\t
T_b}|a_1(t)|\,dt)=O(\ve^2).\tag4.38
$$

By the definition of $u_a^{II}$ in (2.31), one has that
$u_a^{II}(1/\ve)=\ve w_0(1/\ve)$ on $\t\G_{\t\rho_0}^+$ holds
true. Moreover, it follows from Lemma 3.2 that
$$
|Z^\al(u-u_a^{II})|\leq C_{\alpha
b}\varepsilon^{\f{3}{2}}|\ln\ve|(1+t)^{-\f{1}{2}}(1+|t-r|)^{1/2}.\tag
4.39
$$

On the other hand, by [9, Theorem 6.2.1] one has
$$
|\p^{\al}Z^{\beta}\bigl(u_a^{II}(1/\ve)-(\f{\ve}{r^{\f{1}{2}}}F_0)
(r(1/\ve)-1/\ve)\bigr)|\leq C_{\al\beta}\ve^\f{3}{2}.
$$

Therefore,
$$
\align w_1(1/\ve)
&=(r^{\f{1}{2}}\p_t
u)(1/\ve)-\bigl(c(\f{u}{2r^{\f{1}{2}}}+r^{\f{1}{2}}\p_ru)\bigr)(1/\ve)\\
&=(r^{\f{1}{2}}\p_t
u_a^{II})(1/\ve)-\bigl(c(\f{u_a^{II}}{2r^{\f{1}{2}}}
+r^{\f{1}{2}}\p_ru_a^{II})\bigr)(1/\ve)+O(\ve^{3/2})\\ &=\ve
F_0'(r(1/\ve)-1/\ve))(-c(1/\ve)-1)+O(\ve^{3/2}|\ln\ve|)\\ &=-2\ve
F_0'(r(1/\ve)-1/\ve))+O(\ve^{3/2}|\ln\ve|).
\endalign
$$

Note that $|r-t|\leq C+|\t\rho_0|$ and $|u(t,r)|\leq
C_b\ve(1+t)^{-\f{1}{2}}$ on $\t\G_{\t\rho_0}^+$.  Hence,
$$
|r(1/\ve)-1/\ve-\t\rho_0|\leq\int_0^{1/\ve}C_b\ve^2(1+s)^{-1}ds
=C_b\ve^2|\ln\ve|.
$$

Thus one has
$$w_1(\f{1}{\ve})=-2\ve F_0'(\t\rho_0)+O(\ve^{3/2}|\ln\ve|).\tag4.40$$

For later reference, we now provide properties of $u$ when restricted
to $\t\Gamma_{\t\rho_0}^{+}$.

By the definition of $u_a^{II}$ in (2.31), one has, for $t\geq
\ds\f{1}{\ve}$,
$$u_a^{II}(t,x)=\ve\biggl(\chi(\ve t)w_0(t,x)+(1-\chi(\ve
t))r^{-1/2}F_0(\t\rho_0)\biggr)
\quad  \text{on}\enspace \t\Gamma_{\t\rho_0}^{+}.\tag 4.41$$

Moreover, it follows from [9,Lemma 6.2.1] that
$$|\p^{\al}Z^{\beta}(w_0(t,x)-r^{-1/2}F_0(\t\rho_0))|\leq
C_{\al\beta}\ve^{\f{1}{2}}(1+t)^{-1}\quad \text{on}\enspace
\t\Gamma_{\t\rho_0}^{+}.\tag 4.42$$

Substituting (4.42) into (4.41) yields, for $t\geq \ds\f{1}{\ve}$,
$$u_a^{II}=\ve r^{-1/2}F_0(\t\rho_0)+O(\ve^{\f{3}{2}})(1+t)^{-1}\quad
\text{on}\enspace \t\Gamma_{\t\rho_0}^{+}.$$

Together with (4.39), this implies, for $t\geq \ds\f{1}{\ve}$,
$$u=\ve
r^{-1/2}F_0(\rho_0)+O(\ve^{\f{3}{2}})|\ln\ve|(1+t)^{-\f{1}{2}}
+O(\ve^{\f{3}{2}})(1+t)^{-1}\quad
\text{on}\enspace \t\Gamma_{\t\rho_0}^{+}.\tag 4.43$$

\smallskip

Relying on the preparations above, we now prove (4.32).

As $F_0(\t\rho_0)F_0'(\t\rho_0)<0$, without loss of generality we can
assume that $F_0(\t\rho_0)<0$ and $F_0'(\t\rho_0)>0$. One then has,
for $t\geq \ds\f{1}{\ve}$,
$$\text{$w_1(\f{1}{\ve})<0$ \ and \
$a_0(t)=\f{2u+O(u^2)}{2(r(t))^{1/2}c(u)}<0$
on $\t\Gamma_{\rho_0}^{+}$}.$$

Consider the equation for $\t w_1=-w_1$. It follows from (4.36) that
on the characteristic curve $\t\Gamma_{\t\rho_0}^{+}$, for
$t>\ds\f{1}{\ve}$,
$$\f{d \t w_1}{dt}(t,r(t))=-a_0(t)(\t w_1)^2+a_1(t)\t w_1-a_2(t),$$
where $\t w_1(\f{1}{\ve})>0$. By Lemma A.4, one has
$$
\bigl(-\int_{1/\ve}^{\t
T_b}\f{1}{2(r(t))^{\f{1}{2}}c(u)}c'(u)\,dt\bigr)
\exp(-\int_{1/\ve}^{\t T_b}|a_1(t)|\,dt)<(\t w_1(1/\ve)-K)^{-1}.
$$

>From this, together with (4.37)--(4.38), (4.40), (4.43), and
$c'(u)=2u+O(u^2)$, one arrives at
$$
\ve F_0(\t\rho_0)(\ln\t T_b-\ln(1/\ve))(1+O(\ve)) <(-2\ve
F_0'(\t\rho_0)+O(\ve^{3/2}|\ln\ve|))^{-1}\exp(10E\ve^2),
$$
which implies
$$
\overline{\lim_{\ve\rightarrow 0}}\ve^2\ln
T_{\ve}\leq-\f{1}{2F_0(\t\rho_0) F_0'(\t\rho_0)}=\nu_0.
$$

Consequently, (4.32) is shown.

\medskip

{\bf Proof of Theorem 1.1.} Under the assumptions of Theorem 1.1, it
follows from (3.22)--(3.23), (4.28), and (4.32) that the lifespan
$T_{\ve}$ satisfies
$$\ds\lim_{\ve\to 0}\ve \sqrt{T_{\ve}}=\tau_0\quad
\text{when}\enspace c(u)=1+u+O(u^2),$$ and
$$\ds\lim_{\ve\to
0}\ve^2\ln T_{\ve}=\nu_0\quad \text{when}\enspace
c(u)=1+u^2+O(u^3).$$
Thus, we have completed the proof of Theorem 1.1. \qed

\vskip 0.5 true cm \centerline{\bf Appendix. Some useful
lemmas}\vskip 0.5 true cm

{\bf Lemma A.1.} (i) {\it For $\phi(t,r)\in C^{1}$,
$$|\p \phi|\leq \f{2}{1+|t-r|}\sum_{|\beta|=1}|Z^{\beta}\phi|.\tag A.1$$}

(ii) {\it Klainerman fields have the following expressions in
$(\tau,\si)$ coordinates:
$$\cases
\p_t=-\p_{\si}+\ds\f{\ve}{2\sqrt{1+t}}\p_{\tau},\\
\p_{r}=\p_{\si},\\
S=\si\p_{\si}+\ds\f{\ve t}{2\sqrt{1+t}}\p_{\tau},\\
H=-\si\p_{\si}+\ds\f{\ve r}{2\sqrt{1+t}}\p_{\tau}.
\endcases$$}

{\bf Proof.} (i) For $\p_t=\ds\f{tS-rH}{t^2-r^2}$ and
$\p_r=\ds\f{tH-rS}{t^2-r^2}$,
$$(1+|t-r|)(|\p_t\phi|+|\p_r\phi|)\le 2(|S\phi|+|H\phi|
+|\p_t\phi|+|\p_r\phi|),$$
and (A.1) is shown.

(ii) This follows from a direct computation. \qed

\medskip

{\bf Lemma A.2.} {\it If $f(t,x)\in C^1(\Bbb R^{+}\times \Bbb R^{2})$
depends only on $(t, r)$ and $\supp f\subseteq\{(t, x)\colon r\leq
M+t\}$, then
$$
\|(1+|t-r|)^{-1}f\|_{L^2}\leq C\|\p_r f\|_{L^2}.
$$
}

{\bf Proof.} Since $\supp f\subseteq\{r\leq M+t\}$,
$$
f(t, r)=-\int_r^{M+t}\p_r f(t, s)\,ds.
$$
It follows that
$$
\align |f(t, r)|^2 & \leq(\int_r^{M+t}|\p_r f(t,
s)|^2(1+|t-s|)^{1/2}\,ds)
\int_r^{M+t}(1+|t-s|)^{-1/2}\,ds\\
&\leq C(\int_r^{M+t}|\p_r
f(t,s)|^2(1+|t-s|)^{1/2}ds)(1+|t-r|)^{1/2}.
\endalign
$$
Thus,
$$
\align \int_0^{M+t}(1+|t-r|)^{-2}|f(t, r)|^2r\,dr
&\leq C\int_0^{M+t}(\int_r^{M+t}|\p_r f(t,s)|^2(1+|t-s|)^{1/2}\,ds)
  (1+|t-r|)^{-3/2}r\,dr\\
&\leq C\int_0^{M+t}|\p_r f(t,s)|^2(1+|t-s|)^{1/2}\,ds
  \int_0^s(1+|t-r|)^{-3/2}r\,dr\\
&\leq C\int_0^{M+t}|\p_r f(t,s)|^2(1+|t-s|)^{1/2}s\,ds
  \int_0^s(1+|t-r|)^{-3/2}\,dr\\
&\leq C\int_0^{M+t}|\p_r f(t,s)|^2s\,ds,
\endalign
$$
and Lemma A.2 is shown. \qed

\medskip

{\bf Lemma A.3 (Generalized Gronwall's inequality).} {\it
Let $f\in C^1[0, \infty)$, $g, h\in C[0,\infty)$ be non-negative and
$$
\f{df^2(t)}{dt}\leq f(t)g(t)+h(t)f^2(t).
$$
Then
$$
f(t)\leq \biggl(f(0)+\f{1}{2}\int_0^t
g(s)ds\biggr)\exp\bigl(\f{1}{2}\int_0^t h(s)ds\bigr).
$$}

\medskip

{\bf Lemma A.4} ([9, Lemma 1.3.2]). {\it Let w be a solution in $[0,T]$
to the ordinary differential equation
$$
\f{dw}{dt}=a_0(t)w^2+a_1(t)w+a_2(t)
$$
with $a_j$ continuous and $a_0\geq 0$. Let
$$
K=(\int_0^T|a_2(t)|\,dt)\exp(\int_0^T|a_1(t)|\,dt).
$$
Then
$$
(\int_0^Ta_0(t)\,dt)\exp(-\int_0^T|a_1(t)|\,dt)<(w(0)-K)^{-1}
$$
provided that $w(0)>K$.}

\medskip

{\bf Lemma A.5} (Blowup of smooth solution to problem (2.37)). {\it
The smooth solution to\/ {\rm (2.37)} blows up in finite time if
$F_0(\si)\not\equiv 0$.}

{\bf Proof.}  Assume that (2.37) admits a global smooth solution.  Due
to $(F'_0)^2(M)=(F'_0)^2(-\infty)=0$ and $(F'_0)^2\not\equiv 0$, one has
$$F_0'(\si)F_0''(\si)=\bigl(\f12(F_0')(\si)\bigr)'<0\enspace\text{on some
interval $I\subset (-\infty, M)$}.$$
Without loss of generality, we can assume $F'_0(\si)<0$ and
$F''_0(\si)>0$ on $I$.

Let $\Sigma=\{(\tau, \si(\tau, l))\colon \tau\ge 0,\, l\in I\}$, where
$\si(\tau, l)$ stands for the characteristics of (2.37) emanating
from the point $(l,0)$, i.e., $\si(\tau, l)$ satisfies
$$
\cases &\ds\f{d\si(\tau, l)}{d\tau}=G^2(\tau, \si(\tau, l)),\\
&\si(0, l)=l
\endcases\tag A.2
$$

Set $Q(\tau,l)=(\p_{\si}G)(\tau, \si(\tau, l))$ and
$G(\tau,l)=G(\tau, \si(\tau, l))$. It follows from the equation in
(2.37) that
$$Q(\tau, l)=\ds\f{F'_0(l)}{1+F'_0(l)\int_0^{\tau}G(t, l)\,dt}\tag A.3$$
and
$$(\p_{\si}Q)(\tau, \si(\tau, l))\p_l\si(\tau, l)=\ds\f{F''_0(l)-(F'_0)^2(l)
\int_0^{\tau}Q(t, l)\p_l\si(t,l)\,dt} {(1+F'_0(l)\int_0^{\tau}G(t,
l)\,dt)^2}.$$ Therefore,
$$\text{$Q<0$ and $\p_{\si}Q>0$ in $\Sigma$}.\tag A.4$$

Choose $l_i\in I (i=0, 1, 2)$ such that $l_0<l_1<l_2$ and denote
$E_j=\int_{l_j}^{l_{j+1}}(F_0')^2(l)\,dl$ for $j=0,1$. It follows from
conservation of energy for problem (2.37) and (A.4) that, for $j=0,1$,
$$
\align
0<E_j&=\int_{\si(\tau, l_j)}^{\si(\tau, l_{j+1})}Q^2(\tau,
s)\,ds \\
&\le (-Q(\tau, l_j))\int_{\si(\tau, l_j)}^{\si(\tau,
l_{j+1})}(-Q)(\tau, s)\,ds
=-Q(\tau, l_j)\bigl(G(\tau, l_j)-G(\tau,
l_{j+1})\bigr),\endalign$$
which yields
$$G(\tau, l_j)-G(\tau,
l_{j+1})\ge-\ds\f{E_j}{Q(\tau, l_j)},\quad j=0,1.\tag A.5$$

By (A.2), one has that $\si(\tau, l)<M$ holds for $l\in I$ and all
$\tau$. Therefore,
$$\ds\sum_{i=0}^2\int_0^{\infty}G^2(\tau, l_i)\,d\tau
\le 3M+\ds\sum_{i=0}^2|l_i|,$$
which implies that there exists a sequence $\{\tau_k\}\subset [0,
\infty)$ with $\tau_k\to\infty$ as $k\to\infty$ such that
$$G(\tau_k, l_i)\to 0\quad\text{as $k\to\infty$ for $i=0,1,2.$}\tag A.6$$
It then follows from (A.5)--(A.6) that
$$\text{$Q(\tau_k, l_j)\to -\infty$ as $k\to \infty$ for $j=0,1$}.\tag A.7$$

On the other hand, by (A.4), one has
$$\int_0^{\tau_k}G(t, l_0)dt>\int_0^{\tau_k}G(t, l_1)dt.$$
Together with (A.4) and (A.7), this yields as $k\to\infty$
$$-\ds\f{1}{F'_0(l_0)}\ge -\ds\f{1}{F'_0(l_1)}$$
and $F'_0(l_1)\le F'_0(l_0)$. The latter, however,
contradicts the fact that $F'_0(l_0)<F'_0(l_1)$ holds due to
$F''_0(\si)>0$ in $I$ and $l_0<l_1$.

Thus, the proof of Lemma A.5 has been completed. \qed


\bye

\Refs \refstyle{C}

\ref\key 1\by R.K.~Agarwal, D.W.~Halt\paper A modified CUSP scheme
in wave/particle split form for unstructured grid Euler flows \inbook
Frontiers of Computational Fluid Dynamics \eds D.A.Caughey and
M.M.Hafez \yr 1994 \endref

\ref\key 2\by S.~Alinhac \paper Blowup of small data solutions for
a quasilinear wave equation in two space dimensions \jour Ann. of Math.
(2) \vol 149 \pages 97--127 \yr 1999 \endref

\ref\key 3\bysame \paper The null condition for quasilinear
wave equations in two space dimensions I \jour
Invent. Math. \vol 145 \pages 597--618 \yr 2001 \endref

\ref\key 4\bysame \paper An example of blowup at infinity
for quasilinear wave equations \jour Asterisque \vol 284 \pages
1--91 \yr 2003
\endref

\ref\key 5\by D.~Christodoulou\paper Global solutions of nonlinear
hyperbolic equations for small initial data \jour Comm. Pure Appl.
Math. \vol 39 \pages 267--282\yr 1986\endref

\ref\key 6\by R.T.~Glassey; J.K.~Hunter; Zheng, Yuxi \paper
Singularities of a variational wave equation \jour J. Differential
Equations \vol 129 \pages 49--78\yr 1996
\endref

\ref\key 7\by P.~Godin\paper Lifespan of solutions of semilinear
wave equations in two space dimensions \jour Comm. Partial Differential
Equations \vol 18 \pages 895--916\yr 1993 \endref

\ref\key 8\by  L.~H\"{o}mander \paper The lifespan of classical
solutions of nonlinear hyperbolic equations \finalinfo Mittag-Leffler report
No.~5 \yr 1985
\endref

\ref\key 9\bysame \book Lectures on nonlinear
hyperbolic equations \bookinfo Mathematiques \& Applications \vol
26 \publ Springer \publaddr Heidelberg \yr 1997 \endref

\ref\key 10\by A.~Hoshiga\paper The asymptotic behaviour of the
radially symmetric solutions to quasilinear wave equations in two
space dimensions \jour Hokkaido Math. J. \vol 24 \pages 575--615 \yr 1995
\endref

\ref\key 11\bysame \paper The existence of the global
solutions to semilinear wave equations with a class of cubic
nonlinearities in $2$-dimensional space \jour Hokkaido Math. J. \vol 37
\pages 669--688\yr 2008
\endref

\ref\key 12\by J.K.~Hunter; R.A.~Saxton\paper Dynamics  of director
fields \jour SIAM. J. Appl. Math. \vol 51 \pages 1498--1521\yr 1991\endref

\ref\key 13\by F.~John \paper Blow-up of radial solutions of
$u_{tt} = c^2(u_t)\Delta u$ in three space dimensions \jour Mat. Apl.
Comput. \vol 4 \pages 3--18 \yr 1985
\endref

\ref\key 14\by S.Klainerman\paper The null condition and global
existence to nonlinear wave equations \jour Lectures in Applied
Mathematics \vol 23 \pages 193-236 \yr 1986
\endref

\ref\key 15\bysame \paper  Remarks on the global Sobolev
inequalities in the Minkowski space $\Bbb R^{n+1}$ \jour Comm. Pure
Appl. Math. \vol 40 \pages 111-117 \yr 1987
\endref

\ref\key 16\by S.~Klainerman; G.~Ponce\paper Global, small amplitude
solutions to nonlinear evolution equations \jour Comm. Pure Appl. Math.
\vol 36 \pages 133--141\yr 1983\endref

\ref\key 17\by Li, Ta-tsien; Chen, Yun-mei\paper Initial value
problems for nonlinear wave equations \jour Comm. Partial Differential
Equations \vol 13 \pages 383--422\yr 1988 \endref

\ref\key 18\by H.~Lindblad \paper On the lifespan of solutions of
nonlinear wave equations with small initial data \jour Comm. Pure
Appl.  Math. \vol 43 \pages 445--472 \yr 1990 \endref

\ref\key 19\bysame\paper Global solutions of nonlinear
wave equations \jour Comm. Pure Appl. Math. \vol 45 \pages
1063--1096 \yr 1992 \endref

\ref\key 20\bysame \paper Global solutions of quasilinear
wave equations \jour Amer. J. Math. \vol 130 \pages 115--157 \yr
2008 \endref

\ref\key 21\by J.~Metcalfe; C.D.~Sogge\paper Global existence for
high dimensional quasilinear wave equations exterior to
star-shaped obstacles \jour Discrete Contin. Dyn. Syst. \vol 28
\pages 1589--1601\yr 2010 \endref

\ref\key 22\by R.A.Saxton\paper Dynamic instability of the liquid
crystal director \inbook Current progress in hyperbolic systems:
Riemann problems and computations (Brunswick, ME, 1988) \pages
325--330 \bookinfo Contemp.Math. \vol 100 \publ
Amer. Math. Soc. \publaddr Providence, RI \yr 1989\endref

\ref\key 23\by Yin, Huicheng\paper Formation and construction of a
shock wave for 3-D compressible Euler equations with the spherical
initial data \jour Nagoya Math. J. \vol 175 \pages 125--164 \yr
2004 \endref

\ref\key 24\by Yin, Huicheng; Qiu, Qingjiu\paper The blowup of
solutions for 3-D axisymmetric compressible Euler equations \jour
Nagoya Math. J. \vol 154 \pages 157--169 \yr 1999 \endref

\ref\key 25\by Zhang, Ping; Zheng, Yuxi\paper  On the second-order
asymptotic equation of a variational wave equation \jour Proc. Roy.
Soc. Edinburgh Sect. A \vol 132 \pages 483--509 \yr 2002\endref

\ref\key 26\bysame \paper Energy conservative
solutions to a one-dimensional full variational wave system \jour Comm.
Pure Appl. Math. \finalinfo DOI: 10.1002/cpa.20380 (online)\yr 2011\endref

\ref\key 27\by Zheng, Yuxi\paper Existence of solutions to the
transonic pressure-gradient equations of the compressible Euler
equations in elliptic regions \jour Comm. Partial Differential
Equations \vol 22 \pages 1849--1868 \yr 1997\endref

\endRefs
\enddocument